\documentclass[a4paper,11pt]{article}
\usepackage[english]{babel}
\usepackage{theorem}
\usepackage{amssymb}
\usepackage{amsfonts}
\usepackage{amsmath}
\usepackage[all]{xy}
\usepackage{geometry}
\usepackage{ae,aecompl}

\newtheorem{thm}{Theorem}[section]
\newtheorem{prop}[thm]{Proposition}
\newtheorem{lem}[thm]{Lemma}
\newtheorem{cor}[thm]{Corollary}

\makeatletter
\begingroup
\gdef\th@upshape{\normalfont
  \def\@begintheorem##1##2{%
        \item[\hskip\labelsep \theorem@headerfont ##1\ ##2.]}%
\def\@opargbegintheorem##1##2##3{%
   \item[\hskip\labelsep \theorem@headerfont ##1\ ##2\ (##3).]}}
\endgroup
\makeatother

\theoremstyle{upshape}

\newtheorem{eg}[thm]{Example}

\newtheorem{defn}[thm]{Definition}

\title{{\bf Quantum Bundle Description of Quantum Projective Spaces}}
\author{R\'{e}amonn \'{O} Buachalla}
\date{}

\def\a{\alpha}
\def\b{\beta}
\def\g{\gamma}
\def\d{\delta}
\def\e{\varepsilon}

\def\l{\lambda}
\def\Lam{\Lambda}

\def\ta{\theta}

\def\w{\omega}

\def\Om{\Omega}

\def\del{\partial}
\def\DEL{\Delta}

\def\iot{\iota}

\def\bC{{\mathbf C}}
\def\bN{{\mathbf N}}

\def\bZ{{\mathbf Z}}

\def\E{{\cal E}}

\def\Ad{\mathrm{Ad}}

\def\exd{\mathrm{d}}
\def\demo{\noindent \emph{\textbf{Proof}.\ }~}

\def\dim{\mathrm{dim}}

\def\id{\mathrm{id}}

\def\ker{\mathrm{ker}}

\def\proj{\mathrm{proj}}

\def\sgn{\mathrm{sgn}}

\def\spn{\mathrm{span}}

\def\th{\mathrm{th}}

\def\th{^\mathrm{th}}

\def\ver{\mathrm{ver}}

\def\inv{^{-1}}

\def\oby{\otimes}

\def\sseq{\subseteq}
\def\well-defined because
{\widetilde}
\def\wh{\widehat}
\def\ol{\overline}
\def\la{\langle}
\def\\la{\left\langle}
\def\ra{\rangle}
\def\>{\right\rangle}
\def\bs{\backslash}	
\def\mto{\mapsto}

\def\qed{\hfill\ensuremath{\square}\par}
\def\iito{\hookleftarrow}

\def\csun{\bC_q[SU_N]}

\def\cpn{\bC_q[\bC P^{N-1}]}
\def\ccpn{\bC P^{N-1}}
\def\csn{\bC_q[S^{2N-1}]}
\def\ccsn{S^{2N-1}}
\def\usl2{\mathcal{U}_q(\mathfrak{sl}_2)}

\def\wcpn{\Om^1_q(\ccpn)}
\def\1{_{(1)}}
\def\2{_{(2)}}
\def\3{_{(3)}}

\def\alg{algebra~}
\def\algn{algebra}
\def\algs{algebras~}

\def\iff{if, and only if,~}

\def\ncg{noncommutative geometry~}

\def\st{such that~}

\def\wrt{with respect to~}

\DeclareMathOperator{\dt}{det}

\begin{document}

\maketitle

\begin{abstract}
We realise Heckenberger and Kolb's canonical calculus on quantum projective $(N-1)$-space $\cpn$ as the restriction of a distinguished quotient of the standard bicovariant calculus for the quantum special unitary group $\csun$. We introduce a calculus on the quantum sphere $\csn$ in the same way. With respect to these choices of calculi, we present $\cpn$ as the base space of two different quantum principal bundles, one with total space $\bC_q[SU_N]$, and the other with total space $\csn$. We go on to give $\cpn$ the structure of a quantum framed manifold. More specifically, we describe the module of one-forms of Heckenberger and Kolb's calculus as an associated vector bundle to the principal bundle with total space $\csun$. Finally, we construct strong connections for both bundles.
\end{abstract}

\section{Introduction}

The interaction of the theory of quantum groups and their homogeneous spaces with Connes' formulation of noncommutative geometry is a very important, exciting, and active area of contemporary mathematics. The noncommutative geometry of quantum groups is usually discussed in terms of covariant differential calculi (as introduced by Woronowicz in his seminal paper \cite{Wor}). Meanwhile, Connes' operator algebraic formulation of \ncg is most commonly presented in terms of spectral triples \cite{CONN,ENCG}.

One of the major families of quantum group homogeneous spaces is the family of quantum flag manifolds \cite{LR91,Soibel2,TT}. In \cite{HK} Heckenberger and Kolb showed that these spaces admit exactly two finite-dimensional irreducible covariant first-order differential calculi. Classically these correspond to the holomorphic and anti-holomorphic calculi of the manifold. Thus, the calculus obtained by taking their direct sum is a quantum generalisation of the module of complexified differential forms. (Hereafter we refer to this direct sum calculus as the {\em Heckenberger--Kolb calculus}.) This shows that Woronowicz's theory of covariant differential calculi is intimately suited to the study of the geometry of quantum flag manifolds. 

In the spectral triple approach one deals with first-order calculi that arise as subspaces of the \alg of bounded operators on a Hilbert space. The exterior derivative is of the form $\exd f = [D,f]$, for $D$ a self-adjoint operator generalising the classical Dirac operator of a Riemannian spin$^c$ manifold. In \cite{Krah} Kr\"ahmer showed that the Heckenberger--Kolb calculus could be realised in just such a manner. This exciting result has attracted a good deal of attention as it suggests that the spectral triple and covariant calculi approaches to \ncg are compatible.

Unfortunately however, there are several shortcomings of Kr\"ahmer's Dirac operator, not least that its classical limit gives something that is not quite a classical Dirac operator. (See \cite{DDLCP1,NeshTus} for a discussion of some of these shortcomings.) This inspired two papers \cite{DDLCP1,DDCPN} by D'Andrea, D\k{a}browski, and Landi, in which a novel, more explicit, reworking of the Heckenberger--Kolb calculus was developed for the special case of the quantum projective spaces. This reworking made possible the construction of a new Dirac operator realizing the calculus. The operator was built in a manner that modeled the classical picture more closely than Kr\"ahmer's approach, and avoided its failings. Moreover, the construction was shown to satisfy Connes' spectral triple axioms.

In this general context, the aim of the present paper is to offer another reworking of the Heckenberger--Kolb calculus for the quantum projective spaces. This reformulation will be expressed in terms of Brzezi\'nski and Majid's theory of quantum principal bundles \cite{qqgauge,Majleabh,Majprimer}. Classically, every flag manifold can be understood as the base space of a principal bundle. However, thus far the only quantum flag manifold (endowed with the Heckenberger--Kolb calculus) to be presented in quantum principal bundle terms has been the very simplest one: the Podle\'s sphere endowed with the Podle\'s calculus \cite{qqgauge,Maj}. (Recall that the Podle\'s sphere was originally introduced by Podle\'s in \cite{PodlesSphere}, and is one of the best studied examples of a quantum homogeneous space. The Podle\'s calculus was originally introduced in \cite{Podcalc}.) When the general theory was applied to this example, it showed how the calculus naturally decomposes into holomorphic and anti-holomorphic subcalculi. Moreover, it showed how to describe these subcalculi as associated quantum vector bundles, and hence how to describe $\bC_q[\bC P^1]$  as a quantum framed manifold. General methods also demonstrated how to construct a strong connection for the bundle, as well as suggesting a canonical Hodge {$*$-operator}, spin structure, Dirac operator, and Laplace operator. Many of these tools would later be used as basic ingredients in other works. For example, they were employed in the study of the complex geometry of $\bC_q[\bC P^1]$ in \cite{KLVSCP1}, in the construction of guaged Laplacians in \cite{LRZGaugedLaplac,LZ1}, and in the construction of anti-self-dual connections in \cite{LDAntiSD}. We expect that our extension of the quantum principal bundle description of the Podle\'s sphere to include all quantum projective spaces will prove to be of similar use in the study of these spaces. It is anticipated, for instance, that it will allow for a simplification of the approach used in \cite{DDCPN} to construct Dirac operators, and used in \cite{KKCP2,KKCPN} to study complex structures. It also seems likely that our work can be built upon to give a quantum principal bundle description of all quantum flag manifolds endowed with the Heckenberger--Kolb calculus.

\bigskip

The paper is organised as follows: Section $2$ is preliminary. It introduces basic material about differential calculi over Hopf algebras along with the standard bicovariant calculus for a coquasi-triangular Hopf algebra. It also discusses the general theory of quantum principal bundles, quantum framed manifolds, connections, and covariant derivatives for associated quantum vector bundles.

In Section $3$ we recall the basic details of the quantum groups $\bC_q[U_{N}]$ and $\csun$. We then present the coinvariant sub\algs $\cpn = \csun^{U_{N-1}}, \csn = \csun^{SU_{N-1}}$, and $\cpn = \csn^{U_1}$ as Hopf--Galois extensions.

We construct a differential calculus for $\csun$ in Section $4$ by taking a distinguished quotient of the standard bicovariant calculus. In Section $5$ we show how this calculus induces quantum principal bundle structures on our two bundles.

In Section $6$ and $7$ we use the methods of \cite{Maj1,Maj} to frame the calculi on $\csn$ and $\cpn$. We show how the calculus on $\cpn$ naturally decomposes into a direct sum of holomorphic and anti-holomorphic subcalculi. Moreover, we calculate the bimodule relations for the holomorphic and anti-holomorphic subcalculi.

In Section $8$, we construct a strong connection for the bundle $\csun \hookleftarrow \cpn$, and show that it restricts to a strong connection for the bundle \linebreak $\csn \hookleftarrow\cpn$. We show that every associated bundle to  $\csn \hookleftarrow\cpn$ is also an associated bundle to $\csun \hookleftarrow\cpn$, and that the two connections induce the same covariant derivative on such bundles. Finally, we look at the example of the quantum line bundles over $\cpn$ and see how our connections act on them.

\bigskip

The noncommutative complex geometry of the maximal prolongation of $\Om^1_q(\bC P^{N-1})$ will be considered in \cite{MMF2}, following which a K\"ahler structure and a Dirac operator will be introduced.

\section{Preliminaries on Quantum Principal Bundles}

In this section we fix notation and recall the definitions, constructions, and results of the theory of quantum principal bundles that will be used later on. References are provided where proofs or basic details are omitted. 

Let $A$ be an algebra. (In what follows all \algs are assumed to be unital.) A {\em first-order differential calculus} over $A$ is a pair $(\Om^1,\exd)$, where $\Omega^1$ is an $A$-$A$-bimodule and $\exd:A \to \Omega^1$ is a linear map for which holds the {\em Leibniz rule}
\begin{align*}
\exd(ab)=a(\exd b)+(\exd a)b,&  & (a,b,\in A),
\end{align*}
and for which $\Om^1 = \spn_{\bC}\{a\exd b\, | \, a,b \in A\}$. We call an element of $\Om^1$ a $1$-{\em form}.
The {\em universal first-order differential calculus} over $A$ is the pair
$(\Om^1_u(A), \exd_u)$, where $\Om^1_u(A)$ is the kernel of the product map $m: A \otimes A \to A$ endowed
with the obvious bimodule structure, and $\exd_u$ is defined by
\begin{align}\label{univd}
\exd_u: A \to \Omega^1_u(A), & & a \mto 1 \oby a - a \oby 1.
\end{align}
It is not difficult to show that every calculus over $A$ is of the form $(\Omega^1_u(A)/N, \,\proj \circ \exd_u)$, where $N$ is a sub-bimodule of $\Omega^1_u(A)$, and  $\proj:\Omega^1_u(A) \to \Omega^1_u(A)/N$ is the canonical projection. Given two calculi $(\Om,\exd)$ and $(\Om',\exd')$, we say that $(\Om',\exd')$ is a {\em subcalculus} of $(\Om,\exd)$ if there exists a third calculus $(\Om'',\exd'')$ \st $\Om = \Om' \oplus \Om''$ and $\exd = \exd' + \exd''$, where the sum direct sum taken is the obvious one.

Let $H$ be a Hopf \alg with comultiplication $\DEL_H$, counit $\e_H$, antipode $S_H$, unit $1_H$, and multiplication $m_H$ (in what follows we will almost always omit explicit reference to $H$ when denoting these operators). A differential calculus $\Omega^1(A)$ over a left $H$-comodule $A$ is said to be {\em left-covariant} if there exists a left-coaction $\DEL_L:\Omega^1(A) \to H \oby \Omega^1(A)$ \st
\begin{align*}
\DEL_L (a \exd b) = \DEL(a)(\id \oby \exd)\DEL(b), &  & (a,b \in A).
\end{align*}

A calculus over a right $H$-comodule is said to be {\em right-covaraint} if there exists an analogous right-coaction $\DEL_R$. A calculus over a $H$-bicomodule that is both left and right-covariant is said to be {\em bicovariant} if $(\id \oby \DEL_R)\circ \DEL_L = (\DEL_L \oby \id) \circ \DEL_R$. Of course, $H$ itself is a $H$-bicomodule. The left-covariant differential calculi over $H$ were classified in \cite{Wor} as follows: Consider the linear isomorphism
\begin{align*}
s: H \otimes  H \rightarrow H \otimes H,   & &  a \otimes b \mto  aS(b_{(1)}) \otimes b_{(2)},
\end{align*}
with inverse
\begin{align*}
s\inv: H \otimes H \to H \otimes H, & &  a \otimes b \mto ab_{(1)} \otimes b_{(2)}.
\end{align*}
The restriction of $s$ to the universal calculus $\Om^1_u(H)$ gives a linear isomorphism
\begin{align}\label{somegahbyh+}
s: H \otimes H^+  \to \Omega^1_u(H),
\end{align}
where $H^+ =\ker (\e)$ denotes the augmentation ideal of $H$. Now for any right ideal $I_H$ of $H^+$, it can be shown that $s(H \oby I_H)$ is a sub-bimodule of $\Om^1_u(H)$ for which the corresponding calculus $\Om^1(H)$ is left-covariant. Moreover, it can be shown that every left-covariant calculus arises in this way. This correspondence is bijective, meaning that the left-covariant calculi over $H$ are classified by the right ideals of $H^+$. If we denote $\Lambda_H^1 = H^+/I_H$, then it is clear that $s\inv$ descends to an isomorphism between $H \oby \Lambda^1_H$ and $\Omega^1(H)$. In what follows we will usually drop any explicit reference to $s$ and tacitly identify these two spaces. Building upon the classification of left-covariant calculi, it can be shown that bicovariant calculi are in bijective correspondence with the $\Ad_R$-stable right ideals of $H^+$, that is, right ideals $I_H$ \st $\Ad_R(I_H) \sseq I_H \oby H$, where as usual $\Ad_R(h) = h_{(2)} \oby S(h_{(1)})h_{(3)}$, for $h \in H$.

We say that $H$ is {\em coquasi-triangular} if it is equipped with a convolution-invertible linear map $r :H \otimes H\rightarrow \bC$ obeying
\begin{align} \label{coq1}
r(fg\oby h) = r(f\oby h_{(1)})r(g\oby h_{(2)}),            &   & r(f\oby gh) = r(f_{(1)}\oby h)r(f_{(2)}\oby g),
\end{align}
and
\begin{align*} 
 g_{(1)}f_{(1)}  r(f_{(2)}\oby g_{(2)})  =  r(f_{(1)}\oby g_{(1)})f_{(2)} g_{(2)},
\end{align*}
for all $f,g,h \in H$. For any coquasi-triangular Hopf \alg $H$, the {\em quantum Killing form} is the map 
\[
{\mathcal Q}: H \oby H \to \bC ~~~~~~~~~~~~~~~  h \oby g \mto  r(g_{(1)} \oby h_{(1)})r(h_{(2)} \oby g_{(2)}).
\]
If $H$ has a set of generators $\{u^i_j \, |\, i,j = 1, \ldots, N\}$, for some $N \in \bN$, then we can use ${\mathcal Q}$ to define a family of maps $\{Q_{kl}\,|\,k,l = 1, \ldots ,N\}$ by setting 
\[
Q_{kl}: H  \to \bC,  ~~~~~~~~~~~~~~~  h  \mto {\mathcal Q}(h \oby u^k_l).
\]
Using this family of maps, an $N^2$-dimensional representation $Q$ can then be defined by
\begin{align*}
Q: H  \to M_N(\bC)  &  &  h \mto [Q_{kl}(h)]_{kl}.
\end{align*} 
We call $Q$ the {\em quantum Killing representation} of $H$. It can be shown \cite{Maj0} that \linebreak $\ker(Q)^+ = \ker(Q) \cap H^+$ is an $\Ad_R$-stable right ideal of $H^+$, and so, it corresponds to a bicovariant calculus. We call the corresponding calculus the {\em canonical bicovariant calculus} over $H$, and denote it by $\Om_{\text{bc},q}^1(H)$. When $H=\bC_q[SU_2]$, it can be shown that one recovers Woronowicz's $4D_+$ calculus \cite{Wor}. More generally, for $H = \csun$, one recovers the bicovariant calculus introduced by Jur\u{c}o in \cite{Jurc}.

For a right $H$-comodule $V$ with coaction $\DEL_R$, we say that an element $v \in V$ is {\em coinvariant} if $\DEL_R(v) = v \oby 1$, we denote the subspace of all coinvariant elements by $V^H$, and call it the {\em coinvariant subspace} of the coaction. (We define a coinvariant subspace of a left-coaction analogously.) Now for a right $H$-comodule \alg $P$, its coinvariant subspace $M = P^H$ is clearly a sub\alg of $P$. If the mapping
\[
\ver = (m \oby \id) \circ (\id \oby \DEL_R): P \oby_M P \to P \oby H,
\]
is an isomorphism, then we say that $P$ is a {\em Hopf--Galois extension} of $H$.
It is well-known, and not too difficult to show, that this condition is equivalent to exactness of the sequence
\begin{align} \label{qpbexactseq}
0 \longrightarrow P\Om^1_u(M)P {\buildrel \iota \over \longrightarrow} \Om^1_u(P) {\buildrel {ver}\over \longrightarrow} P \oby H^+ \longrightarrow 0,
\end{align}
where $\Om^1_u(M)$ is the restriction of $\Om^1_u(P)$ to $M$, and $\iota$ is the inclusion map (see \cite{Majleabh} for details). Now it is natural to look for a generalisation of this sequence to one using non-universal calculi. This brings us to the central structure used in this paper:
\begin{defn}
A {\em quantum principal $H$-bundle} is a four-tuple $(P,H,N,I_H)$, where $H$ is a Hopf algebra; $P$ a right $H$-comodule \alg \st P is a Hopf--Galois extension of $M=P^H$; $N$  a sub-bimodule of $\Omega^1_u(P)$ determining a right-covariant calculus $\Omega^1(P)$;  $I_H$ an  $\Ad_R$-stable right ideal of $H^+$ determining a bicovariant calculus $\Om^1(H)$; 
for which holds the equality
\begin{align} \label{vercond}
\ver(N) = P \oby I_H.
\end{align}
\end{defn}
We call $P$ the {\em total space}, $H$ the {\em fibre}, and $M$ the {\em base space}. We usually omit explicit reference to the choice of calculi and refer to $(P,H,N,I)$ as the quantum principal $H$-bundle $P \iito M$. It is clear that every Hopf--Galois extension is a quantum principal bundle for the choice of the universal calculus on the total space, and on the fibre. An immediate consequence of the definition is that for any quantum principal bundle $(P,H,N,I_H)$, we have an exact sequence:
\begin{align} 
0 \longrightarrow P\Om^1(M)P {\buildrel \iota \over \longrightarrow} \Om(P) {\buildrel {\ol{ver}}\over \longrightarrow} P \oby \Lambda^1_H \longrightarrow 0,
\end{align}
where $\Om^1(M)$ is the restriction of $\Om^1(P)$ to $M$,  $\iota$ is the inclusion map, and $\ol{ver}$ the descent of $\ver$ to $\Om^1(P)$ (which is well-defined since (\ref{vercond}) holds).

For Hopf \algs $G,H$, a  {\em homogeneous} right $H$-coaction on $G$ is a coaction of the form $(\id \oby \pi) \circ \DEL$, where $\pi: G \to H$ is a surjective Hopf \alg map. We call the coinvariant sub\alg $M=G^H$ of such a coaction a {\em quantum homogeneous space}, and usually denote it by $\pi:G \to H$. Moreover, when $G$ is a Hopf--Galois extension of $M$, we say that $M$ is a {\em Hopf--Galois quantum homogeneous space}. Let us now look at when non-universal choices of calculi give a Hopf--Galois quantum homogeneous space the structure of a quantum principal bundle: The map $s$ can be used to let $\ver$ act on $G \oby G^+$. As is easily seen,
\begin{align}\label{mver}
\ver: G \oby G^+ \mto G \oby H, & & f \oby g = f \oby \pi(g).
\end{align}
Thus, for any left-covariant calculus on $G$ with corresponding right ideal $I_G \sseq G^+$, and left-covariant calculus on $H$ with right ideal $I_H \sseq H^+$, the  requirement (\ref{vercond}) is satisfied \iff $I_H = \pi(I_G)$. Similarly, it is easy to show that $\Omega(G)$ is right-covariant \iff $(\id \oby \pi)(\Ad_R(I_G)) \sseq I_G \oby H$. In this case we have that
\begin{align*}
\Ad_R(\pi(I_G)) & =  (\pi \oby \pi) \Ad_R(I_G)\sseq (\pi \oby \id)(I_G \oby H) = \pi(I_G) \oby H,
\end{align*}
and so, the calculus on $H$ corresponding to $I_H$ is bicovariant. We collect these observations in the following proposition:
\begin{prop}\cite{Maj} \label{comptprop}
Let $\pi:G \to H$ be a Hopf--Galois quantum principal homogeneous space, and $I_G$ a right ideal of $G^+$. If
\begin{align}\label{compt}
(\id \oby \pi) \Ad_R(I_G) \sseq I_G \oby H,
\end{align}
then $(G,H,s(I_G),\pi(I_G))$ is a quantum principal bundle. We call such a quantum principal bundle a {\em quantum principal homogeneous space}.
\end{prop}

An {\em associated bundle} to a quantum principal $H$-bundle $P \hookleftarrow  M$ is a coinvariant left $M$-submodule of the form ${\E = (P \oby V)^H}$, where $V$ is a $H$-comodule and $P \oby V$ is equipped with the tensor product coaction. For a quantum principal
bundle $P \iito M$, a {\em quantum framed manifold} structure for $M$ is a pair $(\E,s)$, where $\E$ is an associated bundle to $P \hookleftarrow  M$, and $s$ is a left $M$-module isomorphism between $\E$ and $\Om^1(M)$ which we call a {\em framing}. For a right $H$-comodule $V$, a {\em soldering form} is a right $H$-comodule map  $\theta: V \to P\Omega^1(M)$ which induces an $M$-module isomorphism
\begin{align*}
s_\theta: (P \oby V)^H \to \Omega^1(M), & &  p \oby v \mto p \theta(v).
\end{align*}
It can be shown that all framings for $\Om^1(M)$ arise in this way from a soldering form \cite{Maj1}. In general, it is not clear how to find a quantum framed manifold structure, or even if one exists. However, in the case of a quantum principal homogeneous space we have the following theorem:
\begin{thm}\label{homfra} \cite{Maj}
For any quantum principal homogeneous space $\pi:G \to H$ with base space $M$, the vector space $V_M=(G^+\cap M)/(I_G\cap M)$ has a well-defined right $H$-comodule structure given by
\begin{align*}
\Delta_M(\ol{v})= \ol{v_{(2)}} \oby S(\pi(v_{(1)})), &  &  (v \in G^+\cap M),
\end{align*}
\wrt which a soldering form is given by
\begin{align} \label{solderingdefn}
\theta(\ol{v})= S(v_{(1)})\exd v_{(2)}.
\end{align}
\end{thm}
In what follows, we will usually denote $M^+= G^+ \cap M$, and $I_{M} = I \cap M$. Moreover, we define the {\em dimension} of $\Om^1(M)$ to be the vector space dimension of $V_M$. 

A {\em connection} for a quatum principal $H$-bundle $P \iito M$ is a left $P$-module projection $\Pi:\Om^1(P) \to \Om^1(P)$ \st $\ker(\Pi) = P\Om^1(M)P$ and 
\[
\DEL_R \circ \Pi = (\Pi \oby \id) \circ \DEL_R.
\]
Connections are in bijective correspondence with linear maps $\w: \Lambda^1_H \to \Om^1(P)$ for which $\ol{\ver} \circ \w = 1 \oby \id$ and $\DEL_R \circ \w = (\w \oby \id) \circ \ol{\Ad_{R,H}}$, where $\ol{\Ad_{R,H}}$ is the descent of $\Ad_{R,H}$ to the quotient $\Lambda^1_H$. We call such a map $\w$ a {\em connection form}. Explicitly, the connection $\Pi_\w$ corrresponding to a connection form $\w$ is given by
\begin{align} \label{pifromom}
\Pi_\w = m  \circ (\id \oby \w) \circ \ver.
\end{align}
For a quantum principal homogeneous space $\pi:G \to H$ connection forms are in turn equivalent to linear maps ${i:\Lambda^1_H \to \Lambda^1_G}$ \st $\ol{\pi} \circ i = \id$ and
\begin{align} \label{bicovcond}
\ol{\Ad_{R,G}} \circ i = (i \oby \id) \circ \Ad_{R,H},
\end{align}
where $\ol{\pi}$ and $\ol{\Ad_{R,G}}$ are defined to be the unique mappings for which the following diagrams are commutative:
\begin{align*}
\xymatrix{ \ker(\e_G) \ar[rr]^{\Ad_{R,G}} \ar[d]_{\proj}  & & \ker(\e_G) \oby G \ar[d]^{\proj \oby \pi}\\
           \Lambda^1_G \ar[rr]_{\ol{\Ad_{R,G}}}                  & & \Lambda^1_G \oby H,}
&  &  \xymatrix{ \ker(\e_G) \ar[rr]^{\pi} \ar[d]_{\proj}  & & \ker(\e_H)  \ar[d]^{\proj}\\
           \Lambda^1_G \ar[rr]_{\ol{\pi}}                  & & \Lambda^1_H.} 
\end{align*}
(Note that $\ol{\Ad_{R,G}}$ is well-defined because (\ref{compt}) is satisfied, while $\ol{\pi}$ is well-defined because $I_H = \pi(I_G)$.) We call such a map $i$ a {\em bicovariant splitting map}. Explicitly, the connection form associated to $i$ is $\w = s \circ i$.  For a more detailed presentation of connections, connection forms, and bicovariant splitting maps see \cite{qqgauge,qmonop,Majleabh,Majprimer}.

A connection $\Pi$ is called {\em strong} if $(\id - \Pi) (\Om^1(P)) \sseq P\Om^1(M)$. Strong connections are important because they allow us to construct covariant derivatives for all the associated bundles of the principal bundle. Recall that if $\E$ is a bimodule over an \alg $A$ and $\Om^1(A)$ is a differential calculus over $A$, then a {\em covariant derivative} for $\E$ is a linear mapping $\nabla: \E \to \E \oby_A \Om^1(A)$ \st
\begin{align*}
 \nabla(sa) = \nabla(s)a + s \oby \exd a, &        &  (s \in \E, a \in A).
\end{align*}
It was shown in \cite{majhaj}, that for any associated bundle $\E$ to a quantum principal bundle $P \hookleftarrow M$, a strong connection $\Pi$ induces a covariant derivative $\nabla$ on $\E$ defined by
\begin{align} \label{covderivative}
 \nabla:\E \to \E \oby_M \Om^1(M), &   & f \mto (\id - \Pi) \exd f,
\end{align}
where we identify $\E \oby_M \Om^1(M)$ with its canonical image in $\Om^1(P)$.

\section{Quantum Spheres and Quantum Projective Spaces as Hopf--Galois Quantum Homogeneous Bundles}

For $q \in (0,1]$ and $\nu : = q-q^{-1}$, let $\bC_q[M_N]$ be the quotient of the free \alg \linebreak ${\bC \la u^i_j,  \,|\, i,j = 1, \ldots , N \ra}$ by the ideal generated by the elements
\begin{align*}
u^i_ku^j_k  - qu^j_ku^i_k, &  & u^k_iu^k_j - qu^k_ju^k_j,                    & &   \; (1 \leq i<j \leq N); \\
  u^i_lu^j_k - u^j_ku^i_l, &  & u^i_ku^j_l - u^i_ku^j_l-\nu u^i_lu^j_k, & &   \; (1 \leq i<j \leq N,\; 1 \leq k < l \leq N).
\end{align*}
These generators can be more compactly presented as 
\begin{align}\label{Rmatrixrels}
\sum_{w,x=1}^N R^{ac}_{wx} u^w_bu^x_d = \sum_{y,z=1}^N R^{yz}_{bd} u^a_yu^c_z, &  & (1 \leq a,b,c,d \leq N),
\end{align}
where, for $H$ the Heaviside step function,
\begin{align}\label{Rmatrix}
{R}^{ik}_{jl} = q^{\d_{ik}}\d_{il}\d_{kj} + \nu H(k-i)\d_{ij}\d_{kl}.
\end{align}
We can put a bi\alg structure on $\bC_q[M_N]$ by introducing a coproduct $\DEL$ and counit $\e$ that act according to $\DEL(u^i_j) =  \sum_{k=1}^N u^i_k \oby u^k_j$, and $\e(u^i_j) = \d_{ij}$. The {\em quantum determinant} of $\bC_q[M_N]$ is the element
\[
\dt_{N} = \sum\nolimits_{\pi \in S_N}(-q)^{\ell(\pi)}u^1_{\pi(1)}u^2_{\pi(2)} \cdots u^N_{\pi(N)},
\]
with  summation taken over all permutations $\pi$ of $N$ elements, and $\ell(\pi)$ the length of $\pi$. As is well-known, $\dt_N$ is a central and grouplike element of the bialgebra. The centrality of $\det_N$ makes it easy to adjoin an inverse $\det_N \inv$. We extend $\DEL$ and $\e$ by setting
$\DEL(\dt_{N} \inv) = \dt_{N} \inv \oby \dt_{N} \inv$, and $\e(\dt_{N} \inv) = 1$, and denote the new bi\alg by $\bC_q[GL_N]$. If we assume that $q$ is real, then we can endow $\bC_q[GL_N]$ with a $*$-\alg structure by defining
\begin{align*}
(\dt_{N} \inv)^* = \dt_{N}, ~~~~ (u^i_j)^* = (-q)^{j-i}\sum\nolimits_{\pi \in S_{N-1}}(-q)^{\ell(\pi)}u^{k_1}_{\pi(l_1)}u^{k_2}_{\pi(l_2)} \cdots u^{k_{N-1}}_{\pi(l_{N-1})},
\end{align*}
where $\{k_1, \ldots ,k_{N-1}\} = \{1, \ldots, N\}\bs \{i\}$ and $\{l_1, \ldots ,l_{N-1}\} = \{1, \ldots, N\}\bs \{j\}$ as ordered sets. Moreover, we can give $\bC_q[GL_N]$ a Hopf $*$-\alg structure by setting ${S(\dt_{N}^{-1}) = \dt_{N}}$, and $S(u^i_j) = \dt_{N}^{-1} (u^j_i)^*$.
We denote this Hopf $*$-\alg by $\bC_q[U_N]$. For $N=1$, we get the Hopf \alg $\bC[U_1]$, where it is usual to denote $u^1_1 = t$, and $\dt_N\inv = t\inv$. If we quotient $\bC_q[U_N]$ by the ideal $\la \dt_{N} - 1 \ra$, then the resulting \alg is again a Hopf  {$*$-\algn}. We denote it by $\bC_q[SU_N]$.
We can put a coquasi-triangular structure $r$ on $\csun$ by defining 
\begin{align}
r(u^i_j \oby u^k_l) = q^{-\frac{1}{N}}R^{ki}_{jl},
\end{align}
where $q^{\frac{1}{N}}$ is some $N\th$-root of $q$. Its convolution-inverse $\ol{r}$ is given by the mapping \linebreak  $\ol{r}(u^i_j \oby u^k_l) = q^{\frac{1}{N}}\ol{R}^{ki}_{jl}$, where 
\begin{align*}
\ol{R}^{ik}_{jl} = q^{-\d_{ik}}\d_{il}\d_{kj} - \nu H(k-i)\d_{ij}\d_{kl}.
\end{align*}
For $N=2$, we get the well-known Hopf algebra $\bC_q[SU_2]$. We usually denote its four generators by  $a = u^1_1,b = u^1_2,c = u^2_1,d = u^2_2$. 

\subsection{Quantum Projective Spaces}

We are now ready to introduce the quantum projective spaces. 
We use a description, introduced in \cite{Mey}, that presents quantum $(N-1)$-projective space as the coinvariant sub\alg of a $\bC_q[U_{N-1}]$-coaction on $\csun$. This sub\alg is a $q$-deformation of the coordinate \alg of the complex manifold $SU_N/ U_{N-1}$. Recall that classically $\bC P^{N-1}$ is isomorphic to $SU_N/ U_{N-1}$.

\begin{defn}
For the surjective Hopf \alg map $\a_N:\bC_q[SU_N] \to \bC_q[U_{N-1}]$ defined by setting $\a_N(u^1_1) = \dt_{N-1} \inv;$ $\a_N(u^1_i)=\a_N(u^i_1)=0$, for $i \neq 1$; and  $\a_N(u^i_j) = u^{i-1}_{j-1},$ for  $i,j=2, \ldots, N$, we have a homogeneous $\bC_q[U_{N-1}]$-coaction on $\bC_q[SU_N]$ given by ${\Delta_{SU_N,\a} = (\id \otimes \a_N) \circ \Delta}$. {\em Quantum projective $(N-1)$-space} $\cpn$ is defined to be the coinvariant subalgebra of $\DEL_{SU_N,\a}$, that is,
\[
\cpn = \{f \in \bC_q[SU_N]\,|\,\DEL_{SU_N,\a}(f) = f \oby 1\}.
\]
\end{defn}

Let us now present this quantum homogeneous space as a quantum principal homogeneous space. We begin by proving a general result:
\begin{lem} \label{hombundlelem}
For quantum homogeneous space $\pi:G \to H$ with base $M$, we have that $G$ is a Hopf--Galois extension of $M$ if $v(1 \oby p) = 0$, for all $p \in \ker(\pi)$, where the map ${v : G \oby G \to G \oby_M G}$ is defined by setting $v(f \oby g) = fS(g_{(1)}) \oby g_{(2)}$, for $f,g \in G$.
\end{lem}

\demo We will establish this result by introducing a map $\ver \inv: G \oby H \to G \oby_M G$ that acts as an inverse for $\ver$ whenever $v(1 \oby p) = 0$, for all $p \in \ker(\pi)$. Let $i:H \to G$ be a linear mapping \st $\pi \circ i = \id$ (such a mapping can always be constructed) and set $\ver \inv = v \circ (\id \oby i)$. We first show that $\ver \circ \ver \inv = \id$: For any $h \in H$,
\begin{align} \label{ververinvsun}
\ver \circ \ver \inv (f \oby h)          & =  \ver (f S(i(h)_{(1)}) \oby i(h)_{(2)})  =  f S(i(h)_{(1)}) i(h)_{(2)} \oby \pi (i(h)_{(3)})\\
                                         & =  f \e(i(h)_{(1)}) \oby \pi (i(h)_{(2)}) =  f \oby \pi(i(h)) =  f \oby h. \notag
\end{align}
We now move on to to showing that $\ver \inv \circ \ver = \id$: For any $x \in G$, the fact that $\pi \circ i = \id$, implies that $i(\pi(x)) = x + p_x$, for some $p_x \in \ker(\pi)$. This means that
\begin{align*}
\ver \inv  \circ \ver (f \oby g) & = \ver \inv (fg_{(1)} \oby \pi(g_{(2)}))  = v(fg_{(1)} \oby i(\pi(g_{(2)}))) \\
                                 & = v(fg_{(1)} \oby g_{(2)}) + v(f g_{(1)} \oby p_{g_{(2)}}) = fg_{(1)}S(g_{(2)}) \oby g_{(3)}\\ 
                                 & = f \e(g_{(1)}) \oby g_{(2)}  = f \oby g.
\end{align*}
We note that $\ver \inv$ does not depend upon our choice for the  map $i$.
\qed

Using this lemma we now give a detailed alternative proof of a result that was originally established in $\cite{Mey}$.

\begin{cor} \label{SUNCPNQPB}
The quantum homogeneous space $\a_N: \csun \to \bC_q[U_{N-1}]$ has a quantum principal bundle structure.
\end{cor}

\demo It is easy to show that any $p \in \ker(\a_N)$ is of the form
\begin{align}\label{alphakergen}
p = \sum_{i=2}^N u^i_1f_i + \sum_{i=2}^N u^1_ig_i, &  & (f_i,g_i, \in \csun).
\end{align}
Now, for any $f \in \csun$, we have
\begin{align*}
v(1 \oby u^i_1f) & = \sum_{k=1}^N S(f_{(1)})S(u^i_k) \oby u^k_1f_{(2)} =  \sum_{k,l=1}^N S(f_{(1)})S(u^i_k) \oby u^k_1S(u^1_l)u^l_1f_{(2)} \\
                 & = \sum_{k,l=1}^N S(f_{(1)})S(u^i_k) u^k_1S(u^1_l) \oby u^l_1f_{(2)} = \sum_{k,l = 1} S(f_{(1)})\e(u^i_1)S(u^1_l) \oby
                 u^l_1f_{(2)} = 0.\\
\end{align*}
Similarily, it can be shown that $v(1 \oby u^1_ig) = 0$, for any $g \in \csun$.
\qed

\subsection{Quantum Spheres}

As mentioned above, the $q$-deformation of $\cpn$ presented in the previous section is based upon the classical isomorphism between $\ccpn$ and $SU_N/U_{N-1}$. The goal of this section is to give an alternative description of $\cpn$ based upon the classical isomorphism between $\cpn$ and $S^{2N-1}/U_1$, where $S^{2N-1}$ is the $(2N-1)$-sphere. We begin by presenting a $q$-deformation of the coordinate \alg of $S^{2N-1}$ which was first introduced in $\cite{SoibelValk}$. This deformation is based upon yet another classical isomorphism, this time the identification of $S^{2N-1}$ and $SU_{N}/SU_{N-1}$.

\begin{defn}
For the surjective Hopf \alg map $\b_N:\bC_q[SU_N] \to \bC_q[SU_{N-1}]$ defined by setting $\b_N(u^1_1) = 1$; $\b_N(u^1_i)=\b_N(u^i_1)=0$, for $i \neq 1$; and $\b_N(u^i_j) = u^{i-1}_{j-1}$, for  $i,j=2, \ldots, N$, we have a homogeneous $\bC_q[SU_{N-1}]$-coaction on $\bC_q[SU_N]$ given by  $\Delta_{SU_N,\b} = (\id \otimes \b_N) \circ \Delta$. The {\em quantum $(2N-1)$-sphere} $\bC_q[S^{2N-1}]$ is the coinvariant subalgebra of $\DEL_{SU_N,\b}$, that is,
\[
\bC_q[S^{2N-1}] = \{f \in \bC_q[SU_N]\,|\,\DEL_{SU_N,\b}(f) = f \oby 1\}.
\]
\end{defn}

The following lemma was originally established in \cite{Mey}. An alternative proof can easily be formulated using the approach of Corollary \ref{SUNCPNQPB} above.

\begin{lem} \label{SUNSNQPB}
The quantum homogeneous space $\b_N: \csun \to \bC_q[SU_{N-1}]$ has a quantum principal bundle structure.
\end{lem}

We would like a description of $\csn$ in terms of generators and relations. We note that, for $i = 1, \ldots, N$, we have
\begin{align*} 
\DEL_{SU_N,\b} (u^i_1)     & = (\id \oby \b_N) (\sum_{k=1}^N u^i_k \oby u^k_1) = \sum_{k=1}^N u^i_k \oby \b_N(u^k_1) = u^i_1 \oby 1,
\end{align*}
and
\begin{align*}
\DEL_{SU_N,\b} (S(u^1_i)) & = (\id \oby \b_N) (\sum_{k=1}^N S(u^k_i) \oby S(u^1_k))  = \sum_{k=1}^N S(u^k_i) \oby \b_N(S(u^1_k)) = S(u^1_i) \oby 1.
\end{align*}
Thus, $u^i_1$ and $S(u^1_i)$ are contained in $\csn$. We will usually denote $u^i_1$ and $S(u^1_i)$ by $z_i$ and $z_i^*$ respectively.
Using representation theoretic methods, it was established in \cite{SoibelValk} that $\csn$ is in fact generated as a \alg by the elements $z_i, \, z_i^*$. It was also shown that a full set of relations is given by
\begin{align}\label{sphererels}
z_iz_j      =  qz_jz_i,   \qquad (i < j);    \qquad  z_i z_j^*   =  q z_j^*z_i,  \qquad (i \neq j); \qquad \qquad\\
 \qquad \qquad z_i z_i^* - z_i^*z_i + q\inv \nu \sum_{k=i+1}^N q^{2(k-i)} z_k z_k^*;   \qquad \sum_{i=1}^N z_i^* z_i=1.  \qquad \qquad \notag
\end{align}
(More easily accessible versions of the proof can be found in \cite{KSLeabh, Chari}.)

We  now introduce a right $\bC[U_1]$-coaction for the quantum $(2N-1)$-sphere and show that $\cpn$ arises as its coinvariant sub\algn. This alternative description of $\cpn$ comes from \cite{Mey}.

\begin{lem} \label{spherelcoaction}
Define a surjective Hopf \alg map $\g_N:\bC_q[SU_N] \to \bC[U_1]$ by setting $\g_N(u^1_1) = t\inv$; $\g_N(u^k_k) = 1$, for $k = 2, \ldots, N-1$; $\g_N(u^N_N) = t$; and  $\g_N(u^i_j) = 0$, for $i,j = 1, \ldots, N$, and $i \neq j$. The map $(\id \oby \g_N) \circ \DEL$ restricts to a $\bC[U_1]$-coaction on $\csn$ which we denote by $\DEL_{\ccsn,\g}$. Moreover, $\bC_q[\bC P^{N-1}]$ is the coinvariant sub\alg of this coaction, that is,
\[
\bC_q[\bC P^{N-1}] = \{f \in \bC_q[S^{2N-1}]\,|\,\DEL_{\ccsn,\g}(f) = f \oby 1\}.
\]
\end{lem}

\demo That we have a $\bC[U_1]$-coaction on $\csn$ is clear from the fact that
\begin{align} \label{ti1det}
\DEL_{\ccsn,\g}(z_i) = (\id \oby \g_N) \sum_{k=1}^N u^i_k \oby u^k_1 = \sum_{k=1}^N u^i_k \oby \g_N(u^k_1) = z_i \oby t^{-1},
\end{align}
and the similarly established $\DEL_{\ccsn,\g}(z^*_i) = z^*_i \oby t$.

Let us now show that $\bC_q[\bC P^{N-1}]$ is the coinvariant sub\alg of $\DEL_{\ccsn,\g}$: For the canonical projection $\d_{N-1}:\bC_q[U_{N-1}] \to \bC_q[SU_{N-1}]$, we have $\d_{N-1} \circ \a_N = \b_N$, and so, the following diagram is commutative:
\begin{displaymath}
\xymatrix{ \bC_q[SU_N]   ~~~   \ar[d]_{\id} \ar[rrr]^{\DEL_{SU_N,\a}~~~~ }     & & & \bC_q[SU_N]    \oby  \bC_q[U_{N-1}] \ar[d]^{\id \oby \d_{N-1}}\\%
           \bC_q[SU_N]   ~~~   \ar[rrr]^{ \DEL_{SU_N,\b}~~~~ }                 & & & \bC_q[SU_N]    \oby  \bC_q[SU_{N-1}].}
\end{displaymath}
It follows that $\cpn$ is contained in $\csn$. Now let $j:\bC[U_1] \to \bC_q[U_{N-1}]$ be the canonical embedding of $\bC[U_1]$ into $\bC_q[U_{N-1}]$ defined by setting $j(t) = \dt_N$ and $j(t \inv) = \dt_N \inv$. Just as in (\ref{ti1det}), it is easy to show that $\DEL_{SU_N,\a}(z_i)=z_i \oby \dt_N \inv$ and  $\DEL_{SU_N,\a}(z_i^*) = z_i^* \oby \dt_N$, and so, we have another commutative diagram:
\begin{align} \label{cpninsphere}
\xymatrix{\bC_q[S^{2N-1}]~~~ \ar[d]_{\iot} \ar[rrr]^{\DEL_{S^{2N-1},\g}~~~~}  & & & \bC_q[S^{2N-1}] \oby \bC[U_1] \ar[d]^{\iot \oby j}\\%
          \bC_q[SU_N]   ~~~                \ar[rrr]^{\DEL_{SU_N,\a}~~~~}      & & & \bC_q[SU_N]    \oby \bC_q[U_{N-1}].}
\end{align}
That $\cpn$ is the coinvariant sub\alg of $\DEL_{\ccsn,\g}$ follows easily from this.

\qed

Now $\DEL_{S^{2N-1},\g}$ induces a $\bZ$-grading for which $\deg(z_i) = -1$ and $\deg(z_i^*) = 1$. Clearly, $\cpn$ arises as the sub\alg of elements of degree zero. Moreover, each homogeneous subspace $\E_p$ of degree $p$ can be understood as an associated vector bundle to $\csn \hookleftarrow \cpn$. To see this, consider the unique $\bC[U_1]$-coaction on $\bC$ for which $\l \mto \l \oby t^{-p}$,  for $\l \in \bC$, and note that it gives a $\bC[U_1]$-coaction on $\csn \oby \bC \simeq \csn$. It is clear that $\E_p$ is the coinvariant subspace of this coaction. Classically, these associated bundles correspond to the line bundles over $\ccpn$, and so, we call them the {\em quantum line bundles}. Finally, we note that dual to the coaction $\DEL_{\ccsn,\g}$, we have a $U_1$-action on $\csn$ defined by
\begin{align*}
e^{i\ta} \triangleright z_i = e^{i\ta} z_i, & & e^{i\ta} \triangleright z_i^{*} = e^{-i\ta} z_i^{*}, & & (\ta \in [0,2\pi)).
\end{align*}
This provides us with an alternative description of $\cpn$ as the invariant sub\alg of a group action. We will, however, not pursue this viewpoint.

Let us now show that $\csn$ is a Hopf--Galois extension of $\cpn$. 
\begin{prop}
With respect to its $\bC[U_1]$-comodule structure, $\csn$ is a Hopf--Galois extension of $\cpn$.
\end{prop}
\demo
Even though $\csn \hookleftarrow \cpn$ is not a quantum homogeneous space, the argument of Lemma \ref{hombundlelem} can still be applied in this case. However, one needs to be careful about the construction of the inverse for $\ver$. Let $i:\bC[U_1] \to \csn$ be the unique linear map for which $i(t^l) = (z_1^*)^l$ and $i(t^{-l}) = z_1^l$, for $l \geq 0$. We then set $\ver \inv = v \circ (\id \oby i)$, where 
\[
v : \csn \oby \csn \to \csun \oby_{\ccpn} \csn
\]
acts as $v(f \oby g) = fS(g_{(1)}) \oby g_{(2)}$, for $f,g \in \csn$. It is routine to check that $\ver \inv(\csn \oby \bC[U_1])$ is contained in $\csn \oby_{\ccpn} \csn$. (The choice of map $i$ is important here because for an arbitrary $i$ the image of $\ver \inv$ is not guaranteed to lie in $\csn \oby_{\ccpn} \csn$.) Using the exact same argument as in Lemma \ref{hombundlelem}, it can now be shown that $\ver \inv$ is inverse to $\ver$ if $v(1 \oby p) = 0$, for all $p \in \ker (\g_N|_{\ccsn})$. That this condition holds is established just as in Corollarly \ref{SUNCPNQPB}.
\qed

\section{The Bicovariant Calculus on $\bC_q[SU_N]$ and a Distinguished Quotient Thereof} \label{Sec:BicCalcQ}

As explained in the preliminaries, for every coquasi-triangular Hopf \alg $H$, there exists a canonical bicovariant differential calculus $\Om^1_{\text{bc},q}(H)$ over $H$, constructed using the quantum Killing representation. In this section we will recall what the calculus looks like in the case of $\csun$; construct a certain quotient of it; and then explain why this quotient is important.

We begin by stating some very useful formulae (given in terms of the coquasi-triangular structure specified in (\ref{Rmatrix})) for the action of $Q$ on certain distinguished elements of $\csun$. We omit the proof which amounts to a routine application of the coquasi-triangular properties given in (\ref{coq1}).
\begin{lem}
For $Q_{kl}$ defined as above, we have the following formulae:
\begin{align}
 Q_{kl}(u^i_j)                & =   ~~~~\, \sum_{a,z=1}^N               q^{-\frac{2}{N}}R^{ik}_{za}R^{za}_{jl}, \label{q1}\\
 Q_{kl}(S(u^g_h))             & =   ~~~~\, \sum_{a,z=1}^N               q^{2(a-h)+\frac{2}{N}}\ol{R}^{ak}_{zh}\ol{R}^{zg}_{al}, \label{qS}\\
 Q_{kl}(u^i_ju^r_s)           & =   ~~~  \sum_{a,b,x,y,z=1}^N         q^{-\frac{4}{N}}R^{rk}_{zb}R^{iz}_{ya}R^{ya}_{jx}R^{xb}_{sl}, \label{q2}\\
 Q_{kl}(u^i_jS(u^g_h))        & =   ~~~  \sum_{a,b,x,y,z=1}^N         q^{2(b-h)}\ol{R}^{bk}_{zh}R^{iz}_{ya}R^{ya}_{jx}\ol{R}^{xg}_{bl}, \label{q1S} \\
 Q_{kl}(u^i_jS(u^g_h)u^r_s)   & =        \sum_{a,b,c,v,w,x,y,z=1}^N   q^{2(b-h)-\frac{2}{N}} R^{rk}_{zc}\ol{R}^{bz}_{yh} R^{iy}_{xa} R^{xa}_{jw} \ol{R}^{wg}_{bv} R^{vc}_{sl}. \label{q1S1}
\end{align}
\end{lem}

We now introduce a basis of $\Lambda^1_{\text{bc},q,SU_N} = \csun^+/\ker(Q)^+$.

\begin{lem}
The set consisting of the elements $\ol{u^1_1 - 1}$;\, $\ol{u^i_j}$,  for $i \neq j$;  and $\ol{u^i_1S(u^1_i)}$,  for $i \neq 1$,
is an $N^2$-dimensional basis of $\Lambda^1_{\text{{\em bc}},q,SU_N}$.
\end{lem}
\demo The lemma is easily proved by using (\ref{q1}) and (\ref{q1S}) to show that, apart from $u^1_1 - 1$, the given representatives of the proposed basis elements are mapped by $Q$ to different elements of the canonical basis of $M_N(\bC)$, while $u^1_1-$ is the only element with non-zero image under $Q_{11}$.
\qed

The next proposition introduces a distinguished right submodule of $\Lambda^1_{\text{bc},q,SU_N}$ which we will use to construct the quotient of $\Om^1_{\text{bc},q}(SU_N)$ mentioned above.

\begin{lem}{\label{decompositionofV}}
The subspace 
\begin{align*}
V_D = \spn_{\bC}\{\,\ol{u^i_j}\,|\, i,j=2, \ldots N;i \neq j\} + \spn_{\bC}\{\,\ol{u^i_1S(u^1_i)}\,|\, i = 2, \ldots N\}
\end{align*}
is a right submodule of $\Lambda^1_{\text{bc},q,SU_N}$. Moreover, in the quotient $\Lambda^1_{\text{bc},q,SU_N}/V_D$, the subspaces 
\begin{align*}
V_+=\spn_{\bC}\{\ol{u^i_1}\,|\,i=2,\ldots , N\}, & & V_-= \spn_{\bC}\{\ol{u^1_i}\,|\,i=2,\ldots , N\},
\end{align*}
are right submodules.
\end{lem}

\demo That $V_D$ is a right submodule is equivalent to the statement
\begin{eqnarray*}
Q_{1l}(u^i_ju^r_s) = Q_{k1}(u^i_ju^r_s) = Q_{1l}(u^i_1S(u^1_i)u^r_s) = Q_{k1}(u^i_1S(u^1_i)u^r_s) = 0,
\end{eqnarray*}
for all $1 \leq k,l,r,s \leq N$. This is easily proved using (\ref{q2}) and (\ref{q1S1}). That $V_+$ and $V_-$ are right submodules is equivalent to the statement 
\[
Q_{k1}(u^i_1u^r_s) = Q_{1l}(u^1_ju^r_s) = 0,
\]
for all $2 \leq i,j \leq N,$ and $1 \leq k,l \leq N$. This is easily proved using (\ref{q2}). \qed

\begin{cor}
We have a well-defined quotient calculus
\begin{eqnarray}\label{eqn:quotientcalcsln}
\Om^1_q(SU_N)   & = \Om^1_{\text{bc},q}(SU_N)/s(\csun \oby V_D).
\end{eqnarray}
Its corresponding ideal is $I_{SU_N} = \ker(Q)^+ + D_1 + D_2$, where 
\begin{align*}
D_1 = \spn_{\bC}\{u^i_1S(u^1_i)\,|\, i = 2, \ldots, N\}, & \text{~\, and } & D_2= \spn_{\bC}\{u^i_j \,|\, i,j=2 \ldots, N ; i \neq j\}.
\end{align*}
Moreover, for $i=1, \ldots N-1$, the elements 
\begin{align*}
e_{i}^- = \ol{u^{i+1}_1},    &     &  e^0 = \ol{u^1_1 - 1},   &   & e_{i}^+ =  \ol{u_{i+1}^1},
\end{align*} 
form a $(2N-1)$-dimensional left-module basis of $\Om^1_q(SU_N)$.
\end{cor}

We recognise that the dimension of the differential calculus we have chosen for $\csun$ is significantly less than the classical value, for $N>2$. However, we are not interested in $\Om^1_q(SU_N)$ as a quantum deformation in itself. Instead, we will view it as a useful mathematical tool to be exploited in our efforts to investigate the geometry of $\csn$ and $\cpn$. As we shall show below, the calculi that $\Om^1_q(SU_N)$ induces by restriction on $\bC_q[S^{2N-1}]$ and $\cpn$ both have classical dimension; we denote them by $\Om^1_q(S^{2N-1})$ and $\Om^1_q(\bC P^{N-1})$ respectively. By contrast, if we were to consider the calculi that $\Om^1_{\text{bc},q}(SU_N)$ restricts to on $\cpn$ and $\csn$, we would find that they have dimensions $N^2$ and $N^2-1$ respectively, values which are much larger than the classical ones. (The three-dimensional calculus induced by $\Omega^1_{\text{bc},q}(SU_2)$ on $\bC_q[\bC P^1]$ was thoroughly investigated in \cite{Brainlandi}.)

\begin{eg}
Let us look now at the case of $N=2$: The ideal $I_{SU_2}$ corresponding to $\Om^1_q(SU_2)$ is generated by the six elements
\begin{align*}
 (a-q)(a-1), ~~bc,  ~~  b^2, ~~ c^2, ~~ (a-q)b, ~~ (a-q)c,
\end{align*}
or equivalently by the six elements 
\begin{align*}
a+qd-(q+1), ~~bc,  ~~  b^2, ~~ c^2, ~~ (a-q)b, ~~ (a-q)c,
\end{align*}
A three-dimensional basis of left-invariant forms is given by
\begin{align*}
e^+_1 = a\exd c -qc\exd a, &  & e^0 = d\exd a - q\inv b \exd c,   & & e^-_1 = d \exd b - q \inv b \exd d.
\end{align*}
While the exterior derivative acts according to
\begin{align*}
\exd a = a e^0 + b e^+_1, & &      \exd b = a e^-_1 - q\inv be^0, & &   \exd c = ce^0 + d e^+_1, & &    \exd d =  c e^-_1 - q \inv de^0.
\end{align*}
Finally, in matrix form, the right module relations are given by:
\begin{equation*}
e^0\left(
\begin{array}{cc}
a & b\\
c & d
\end{array} \right) = \left( \begin{array}{cc}
                     qa & q \inv b\\
                     qc & q \inv d
                     \end{array} \right)e^{0} + (q-1)\left( \begin{array}{cc}
                     b & 0\\
                     d & 0
                     \end{array} \right) e^{+} + (q - 1)\left( \begin{array}{cc}
                     0 &  a\\
                     0 &  c
                     \end{array} \right)e^{-},
\end{equation*}
\begin{equation*}
e^{\pm}\left(
\begin{array}{cc}
a & b\\
c & d
\end{array} \right) = \left( \begin{array}{cc}
                     a & b\\
                     c & d
                     \end{array} \right)e^{\pm}.
\end{equation*}
Since $\Om^1_q(SU_2)$ is a three-dimensional calculus, it is natural to ask whether or not it is isomorphic to Woronowicz's well-known $3D$ calculus \cite{Wor}. Recall that the ideal corresponding to the the $3D$ calculus is generated by the elements 
\begin{eqnarray*}
a+q^{-2}d-(1+q^{-2}),~ bc, ~ b^2,~c^2,~(a-1)b,~(a-1)c.
\end{eqnarray*}
Using (\ref{q1}) and (\ref{q2}), it is easy to show that $\ol{a+q^{-2}d-(1+q^{-2})}$, \,$\ol{(a-1)b}$ \, and \,$\ol{(a-1)c}$ are all non-zero elements of $\Lambda^1_{SU_2}$. Thus, the two calculi cannot be isomorphic. Alternatively, one can observe that since $(a-q)b - (a-1)b = (1-q)b$, any ideal containing both $(a-q)b$ and $(a-1)b$ will also contain $b$. Since this is not the case for either ideal, they cannot be equal. Moreover, a similar argument will show that $\Om^1_q(SU_2)$ is not isomorphic to any of the {three-dimensional} calculi presented in \cite{SSleftCovSL2}.
\end{eg}

\section{Quantum Principal Bundles Structures Induced by the Calculus $\Om^1_q[SU_N]$}

In this section we will see how the calculus on $\csun$ introduced above induces quantum principal bundle structures on the Hopf--Galois extensions $\a_N:\csun \to \cpn$,\, $\b_N:\csun \to \csn$, and  $\csn \hookleftarrow \cpn$. We will then produce explicit descriptions of the calculi on the fibres $\bC_q[U_{N-1}]$ and $\bC_q[U_{1}]$.

\begin{prop} \label{prop:suntcpn}
It holds that $(\csun,\bC_q[U_{N-1}],I_{SU_N},\a_N(I_{SU_N}))$ is a quantum principal homogeneous space.
\end{prop}

\demo We have already proved that $\a_N:\csun \to \cpn$ is a Hopf--Galois principal homogeneous space. Thus, Proposition \ref{comptprop} tells us that all we need to show is that (\ref{compt}) holds for $I_{SU_N}$. Recall that $I_{SU_N} = \ker(Q)^+ + D_1 + D_2$, where $D_1=\spn_{\bC}\{u^i_1S(u^1_i)\,|\, i = 2, \ldots , N\}$, and $D_2 = \spn_{\bC}\{u^i_j \,|\, i \neq j; i,j =2, \ldots , N\}$. Now since $\ker(Q)^+$ is an $\Ad_R$-stable ideal, it is clear that 
\[
(\id \oby \a_N)\Ad_R(\ker(Q)^+) \sseq \ker(Q)^+ \oby \bC_q[U_{N-1}].
\]
For $D_1$, we begin by noting that
\begin{align*}
(\id \oby \a_N)\Ad_R(u^i_1S(u^1_i)) & = \sum_{a,b,c,d=1}^N u^a_bS(u^c_d) \oby \a_N(S(u^i_aS(u^d_i))u^b_1S(u^1_c))\\
                                    & = \,\, \, \sum_{a,d=2}^N ~~ u^a_1S(u^1_d) \oby S(u^{i-1}_{a-1}S(u^{d-1}_{i-1})) \dt_{N-1} \inv \dt_{N-1}\\
                                    & = \,\, \, \sum_{a,d=2}^N ~~ u^a_1S(u^1_d) \oby S(u^{i-1}_{a-1}S(u^{d-1}_{i-1})).
\end{align*}
For $a=d$, we have $u^a_1S(u^1_a) \in D_1$ by definition. For $a \neq d$, it is easy to use (\ref{q1S}) to show that $Q(u^a_1S(u^1_d))$ is equal to a scalar multiple of $Q(u^a_d)$, and so, $u^a_1S(u^1_d) \in \ker (Q)^+ + D_2$. This means that we have $(\id \oby \a_N)\Ad_R(D_1) \sseq I_{SU_N} \oby \bC_q[U_{N-1}]$. Turning now to $u^i_j \in D_2$, we see that
\[
(\id \oby \a_N)\Ad_R(u^i_j) = (\id \oby \a_N)(\sum_{k,l=1}^N u^k_l \oby S(u^i_k)u^l_j) = \sum_{k,l=2}^N u^k_l \oby S(u^{i-1}_{k-1})u^{l-1}_{j-1}.
\]
For $k \neq l$, we have $u^k_l \in D_2$ by definition. It remains to show that  $\sum_{k=2}^N u^k_k \oby S(u^{i-1}_{k-1})u^{k-1}_{j-1}$ is contained in $I_{SU_N} \oby \bC_q[U_{N-1}]$. Let $\{e_p\,| \,p = 1, \ldots, M \}$ be a basis of the subspace of $\bC_q[U_{N-1}]$ spanned by the elements $\{S(u^{i-1}_{k-1})u^{k-1}_{j-1}) \,|\, k = 2, \ldots ,N\}$, and let $\l_{pk}$ be the unique constants for which
\begin{align*}
\sum_{k=2}^N u^k_k \oby S(u^{i-1}_{k-1})u^{k-1}_{j-1} =  \sum_{k=2}^N ( u^k_k \oby (\sum_{p =1}^M \l_{pk}e_p)) = \sum_{k=2}^N \sum_{p =1}^M \l_{pk} u^k_k \oby e_p.
\end{align*}
We see that $\sum_{k=2}^N u^k_k \oby S(u^{i-1}_{k-1})u^{k-1}_{j-1}$ is contained in $I_{SU_N} \oby \bC_q[U_{N-1}]$ \iff $\sum_{k=2}^N \l_{pk} u^k_k$ is contained in $I_{SU_N}$, for $p = 1, \ldots ,M$. This is equivalent to requiring that
\begin{align*}
\sum_{k=2}^N \l_{pk} Q_{1l}(u^k_k) = \sum_{k=2}^N \l_{pk} Q_{l1}(u^k_k) = \sum_{k=2}^N \l_{pk} Q_{11}(u^k_k) = 0, &  &  (l = 2, \ldots , N).
\end{align*}
Using (\ref{q1}) it is easy to show that $Q_{l1}(u^k_k) = Q_{1l}(u^k_k) = 0$, for all $l = 2,\ldots, N$, and that $Q_{11}(u^k_k) = q^{-\frac{2}{N}}$, for all $k \geq 2$. Thus, $\sum_{k=2}^N \l_{pk} u^k_k \in I_{SU_N}$ \iff $\sum_{k=2}^N \l_{pk} = 0$. That this is true follows from
\begin{align}\label{lambdasum}
\sum_{k=2}^N \l_{pk} & = \sum_{k=2}^N \wh{e}_p(S(u^{i-1}_{k-1})u^{k-1}_{j-1}))  = \wh{e}_p(\e(u^{i-1}_{j-1}))= \e(u^i_j) = 0,
\end{align}
where $\{\wh{e}_p\,| \,p = 1, \ldots, M \}$ denotes the dual basis of  $\{e_p\,| \,p = 1, \ldots, M \}$. (Note that the summation in (\ref{lambdasum}) takes place in $\bC_q[U_{N-1}]$.)
\qed

We now introduce a basis of $V_{\ccpn}$, and hence show that $\Om^1_q(\ccpn)$ has classical dimension:
\begin{lem} \label{lem:cpnbasis}
The set $\{\ol{z_{i1}},\, \ol{z_{1i}} \,|\, i = 2, \ldots,N\}$ is a basis of $V_{\ccpn}$.
\end{lem}
\demo 
Let us identify $V_{\ccpn}$ with its canonical image in $\Lambda^1_{SU_N}$. First we note that $\cpn^+$ is generated as an ideal of $\cpn$ by the elements $z_{11}-1$, and $z_{ij}$, for $(i,j) \neq (1,1)$. Using (\ref{q1S}), it is easy to show that
\begin{align}\label{ztou}
\ol{z_{i1}}=q^{\frac{2}{N}-1}e^+_{i-1}, & & \ol{z_{1i}} = -q^{3+\frac{2}{N}-2i}e^-_{i-1}.
\end{align}
Thus, $e^+_{i-1}$ and $e^-_{i-1}$ are contained in $V_{\bC P^{N-1}}$, for all $i = 2, \ldots, N$. It remains to show that $e^0$ is not contained in $V_{\bC P^{N-1}}$: Formula (\ref{q1S}) can be used to show that $\ol{z_{11} - 1} = 0$, and  $\ol{z_{ij}} = 0$, for all  $i,j \geq 2$. The fact that $V_-$ and $V_+$ are right submodules then implies that $e^0$ is not an element of $V_{\bC P^{N-1}}$. \qed

We are now ready to describe the calculus on the fibre $\bC_q[U_{N-1}]$:
\begin{lem} \label{un1projcalc}
The calculus $\Om^1_q(U_{N-1})$ is one-dimensional with generator $\exd(\dt_{N-1})$ and relations
\begin{align*}
\exd(\dt_{N-1})u^i_j            & =   q^{-\frac{2}{N}} u^i_j\exd (\dt_{N-1}),\\
\exd (\dt_{N-1})\dt_{N-1} \inv  &  =  q^{2-\frac{2}{N}}\dt_{N-1}\inv \exd(\dt_{N-1}). 
\end{align*}
\end{lem}

\demo
Recall that $\Om^1_q(SU_N)$ is isomorphic to ${\csun \oby (V_+ \oplus \bC e^0 \oplus V_-)}$. From (\ref{mver}) it is clear that $ \ker(\ver) = \csun \oby (V_+ \oplus V_-)$ . Thus, since we are dealing with a quantum principal bundle, $\ver$ must map $\csun \oby (\bC e^0)$ isomorphically to $\csun \oby \Lambda^1_{U_{N-1}}$. This means that $\Lam^1_{U_{N-1}}$ is spanned by $\ol{\a_N(u^1_1-1)} = \ol{\det_{N-1} \inv -1}$, or alternatively by 
\[
\ol{\dt_{N-1}-1}(-\dt_{N-1}) = \ol{\dt_{N-1} -1}.
\]
Thus, $\exd(\dt_{N-1}) = \dt_{N-1} \oby (\ol{\dt_{N-1} -1})$  generates the calculus as a left-module over $\csun$.

The commutation relations with the generator $\exd(\dt_{N-1})$ are established just as for the simpler relations established in Lemma \ref{lem:u1projcalc} below.
\qed

\bigskip

It can also be shown that $\Om^1(SU_N)$ induces a quantum principal bundle structure on the Hopf--Galois extension $\csun \hookleftarrow \csn$. Since the proof is identical to the proof of Proposition \ref{prop:suntcpn}, we will simply state the result:
\begin{prop}
It holds that $(\csun,\bC_q[SU_{N-1}],I_{SU_N},\b_N(I_{SU_N}))$ is a quantum principal homogeneous space.
\end{prop}

We now introduce a basis of $V_{\ccsn}$, and hence show that $\Om^1_q(\ccsn)$ has classical dimension:
\begin{lem} \label{lem:spherebasis}
The set $\{\ol{z_{1}-1},\,\,\,\ol{z_i},\,\,\, \ol{z_i^*}) \,|\, i = 2, \ldots,N\}$ is a basis for $V_{\ccsn}$.
\end{lem}
\demo 
Identifying $V_{\ccsn}$ with its canonical image in $\Lambda^1_{SU_N}$, we have that $e^+_{i-1}=\ol{z_i}$, and $e^0 = \ol{z_1-1}$, are contained in $V_{\ccsn}$. Using (\ref{qS}) it is easy to show that 
\begin{align}\label{zstartou}
\ol{z_i^*} = {\ol{S(u^1_i)} = -q^{1+\frac{4}{N}-2i} \ol{u^1_i}},
\end{align}
and so, $e^-_{i-1}$ is also contained in $V_{S^{2N-1}}$. It follows that our proposed basis forms a basis for $\Lambda^1_{SU_N}$, and a fortiori a basis for $V_{\ccsn}$.
\qed
Using an argument analogous to that found in Lemma \ref{un1projcalc}, it is now easy to show that the calculus on $\bC_q[SU_{N-1}]$ corresponding to the ideal $\b_N(I_{SU_N})$ is trivial. (We note that Theorem \ref{homfra} still holds in this case, as a careful reading of the original proof will verfiy.)

\bigskip

Finally, we show that $\Om^1_q[SU_N]$ induces a quantum principal structure on \linebreak $\csn \hookleftarrow \cpn$:

\begin{prop}
It holds that
\[
(\csn,\, \cpn, \,I_{\ccsn}, \,\g_N(I_{\ccsn}))
\]
is a quantum principal bundle.
\end{prop}
\demo
Since we have already shown that $\csn$ is a Hopf--Galois extension of $\cpn$, it suffices to show that $\ver(N_{\ccsn}) = \csn \oby I_{U_1}$; that $\Omega^1_q(\ccsn)$ is right-covariant; and that $I_{U_1}$ is $\Ad_R$-invariant. We begin with the requirement on $\ver(N_{\ccsn})$: It is easily seen from the relations in (\ref{sphererels}) that $\csn$ is equal to the sum of the two subspaces $A = \bC[z_1] + \bC[z^*_1]$ and $B = \bC\la z_i,z^*_i | i = 2, \ldots ,N\ra$. Thus, for any any $\ol{v} \in I_{\ccsn}$, we have a decomposition $\ol{v} = \ol{v_A} + \ol{v_B}$, where $v_A \in A^+, v_B \in B$. Now if we had $v_A,v_B \notin I_{\ccsn}$, then the images of $A^+$ and $B$ in $V_{\ccsn}$, would have non-trivial intersection. However, this is not the case: The fact that $V_+$ and $V_-$ are right submodules of $V_{\ccsn}$ directly implies that 
$\{ \ol{b} \,|\, b \in B \} = V_+ \oplus V_-$, while formulae (\ref{q2}), and (\ref{q1S}), can easily be used to show that $\{\ol{a} \,|\, a \in A^+\} = \bC e^0$. Hence, we have $I_{\ccsn} = I_A \oplus I_B$, where $I_A = I_{\ccsn} \cap A$, and $I_B = I_{\ccsn} \cap B$.
Routine calculation will show that every element of $\csn \oby A^+$ is $\bC_q[SU_{N-1}]$-coinvariant, and so,
\begin{align*}
N_{\ccsn} & = (\csun \oby I_{\ccsn})^{SU_{N-1}} = (\csun \oby (I_{A}\oplus I_B))^{SU_{N-1}}\\
          & = (\csun \oby I_{A}) \oplus (\csun \oby I_B)^{SU_{N-1}}.
\end{align*}
Just as for the general homogeneous case, it is easy to show that $\ver$ acts on \linebreak $(\csun \oby \csn^+)^{SU_{N -1}}$ as $\id \oby \g_N$. Thus,
\begin{align*}
\ver(N_{\ccsn}) & = (\id \oby \g_N)(\csun \oby I_A) + (\id \oby \g_N)(\csun \oby I_B)^{SU_{N-1}}\\
                & = \csun \oby \g_N(I_A) = \csun \oby \g_N(I_{\ccsn}),
\end{align*}
as required.

Let us now move on to establishing the right covariance of $\Om^1_q(\ccsn)$. It is easy to see that there exists a unique \alg map $\zeta_{N-1}:\bC_q[U_{N-1}] \to \bC[U_1]$ for which $\zeta_{N-1} \circ \a = \g_N$. Now
\begin{align}\label{gammainvofomsun}
(\id \oby \g_N) \circ \Ad_R(I_{SU_N}) &  = (\id \oby \zeta_{N-1}) \circ (\id \oby \a_N) \circ \Ad_R(I_{SU_N})\\
                                      &  \sseq (\id \oby \zeta_{N-1})(I_{SU_N} \oby \bC_q[U_{N-1}]) \notag\\
                                      &  = I_{SU_N} \oby \bC_q[U_{1}].\notag
\end{align}
Thus, we must have that $\DEL_{\ccsn,\g}(N_{SU_N}) \sseq N_{SU_N} \oby \bC[U_1]$. Since it is clear that $\DEL_{\ccsn,\g}(\Om^1_u(\ccsn))$ is contained in  $\Om^1_u(\ccsn) \oby \bC[U_1]$, and since by construction $N_{\ccsn} = N_{SU_N} \cap \Om^1_u(\ccsn)$, this means that $\DEL_{\ccsn,\g}(N_{\ccsn}) \sseq N_{\ccsn} \oby \bC[U_1]$, and so, $\Om^1_q(\ccsn)$ is right-covariant.

It is routine to show that every element of $\csn$ is coinvariant under the action $(\id \oby \g_N) \circ \Ad_R$. Together with (\ref{gammainvofomsun}) this gives us that 
\begin{align*}
(\id \oby \g_N) \circ \Ad_R(I_{\ccsn}) \sseq I_{\ccsn} \oby \bC[U_1].
\end{align*}
With this result the bicovariance of $\Om^1(U_1)$ can now be established just as in the general homogeneous case.
\qed

Let us now describe the calculus $\Om^1(U_1)$:
\begin{lem} \label{lem:u1projcalc}
For the bundle $\csn \hookleftarrow \cpn$, the calculus $\Om^1_q(U_1)$ is one-dimensional with $\exd t$ as a generator. Moreover, we have the relation
\begin{align}\label{u1calcrels}
\exd t.t = q^{\frac{2}{N}-2} t\exd t.
\end{align}
\end{lem}
\demo Using methods exactly analogous to those of Lemma \ref{un1projcalc}, it can be shown that  $\Lambda^1_{U_1}$ is a one-dimensional space spanned by the element $\ol{t  -1}$. Since $\exd t = t \oby \ol{t-1}$, it must generate $\Om^1_q(U_1)$ as a left-module over $\bC[U_1]$. To establish (\ref{u1calcrels}) we first note that $(\exd t)t =  t^2 \oby \ol{(t^2-t)}$. Since $\Lambda^1_{U_1}$ is one-dimensional, we must have $\ol{t^2 - t} = \l\ol{t - 1}$, for some $\l \in \bC$. With a view to finding $\l$, we define $\ol{\g_N}$ to be the unique map for which the following diagram is commutative:
\begin{displaymath}
\xymatrix{&\csn^+              \ar[rrr]^{\g_N} \ar[d]_{\proj_{S^{2N-1}}} & & & \bC[U_1]^+ \ar[d]^{\proj_{U_1}} \\
         &V_{S^{2N-1}}         \ar[rrr]^{\ol{\g_N}}                                      & & & \Lambda^1_{U_1}.}
\end{displaymath}
We see that
\begin{align*}
\ol{t^2-t}   = \ol{\g_N(S(u^1_1)^2 - S(u^1_1))}  = \ol{\g_N}(\ol{S(u^1_1)^2-S(u^1_1)}).
\end{align*}
By evaluating the action of $Q$ on $S(u^1_1)^2-S(u^1_1)$, one can show that
\[
\ol{S(u^1_1)^2-S(u^1_1)} = q^{\frac{2}{N}-2}(\ol{S(u^1_1)-1}).
\]
Thus, since $\ol{\g_N}(\ol{S(u^1_1)-1}) = \ol{t-1}$, we must have that $\l = q^{\frac{2}{N}-2}$. Relation (\ref{u1calcrels}) now follows from
\[
(\exd t)t = t^2 \oby \ol{(t^2-t)} = q^{\frac{2}{N}-2} t^2 \oby \ol{t-1} = q^{\frac{2}{N}-2}t\exd t.
\]

Alternatively, the above method can be used to show that $(\exd t\inv)t = q^{\frac{2}{N}-2}t\exd t \inv$, after which (\ref{u1calcrels}) can be concluded from the relation $\exd t = t(\exd t\inv)t$. This approach has the advantage of being computationally simpler. \qed

\section{A Framing for the Odd-Dimensional Quantum Spheres} \label{Sec:Sphere}

In this section we will look at the framing of the calculus $\Om^1_q(S^{2N-1})$ given by Theorem \ref{homfra}. Let us begin by calculating the action of the soldering form on the basis elements of $V_{S^{2N-1}}$ given in Lemma \ref{lem:spherebasis}:
\begin{align*}
\ta(\ol{z_i}) = \ta(\ol{u^i_1})    = \sum_{k=1}^N S(u^i_k)\exd u^k_1 = \sum_{k=1}^N S(u^i_k)\exd z_k,
\end{align*}
\begin{align*}
\ta(\ol{z_1-1}) = \ta(\ol{u^1_1-1})  = \sum_{k=1}^N S(u^1_k)\exd u^k_1 = \sum_{k=1}^N S(u^1_k)\exd z_k,
\end{align*}
\begin{align*}
\ta(\ol{z_i^*}) = \ta(\ol{S(u^1_i)}) = \sum_{k=1}^N S^2(u^k_i)\exd (S(u^1_k)) = \sum_{k=1}^N q^{2(k-
i)}u^k_i\exd z^*_k.
\end{align*}
Moreover, as noted in Lemma \ref{lem:spherebasis}, we have $\ta(\ol{z_i})=e^+_i, ~ \ta(\ol{z_1-1}) = e^0$, and \linebreak $\ta(\ol{z_i^*}) = -q^{1+\frac{4}{N}-2i}e^-_i$.

We now present a decomposition of $\Om^1_q(S^{2N-1})$ into a direct sum of vector spaces. While we will not use this decomposition in our later work, it is of interest in itself as a quantum generalisation of a classical fact.
\begin{lem} \label{decompositionofOm}
Denoting
\begin{align*}
\Om^1_+(S^{2N-1}) = (\csun \oby V_+)^{SU_{N-1}}, &  & \Om^1_0(S^{2N-1}) = (\csun \oby (\bC e^0))^{SU_{N-1}},\\ 
\Om^1_-(S^{2N-1}) = (\csun \oby V_-)^{SU_{N-1}},  & 
\end{align*}
we have the direct sum decomposition of vector spaces
\[
 \Om_q^1(S^{2N-1}) = \Om^1_-(S^{2N-1}) \oplus \Om^1_0(S^{2N-1}) \oplus \Om^1_+(S^{2N-1}).
\]
\end{lem}

\demo For any $\w \in \csn \oby V_{S^{2N-1}}$, we have the decomposition $\w = w_+ + \w_0 + \w_-$, with $\w_0 \in \bC e^0$, and $w_i \in \csun \oby V_i$,\, for $i = +,-$.  Denoting the $\bC_q[SU_N]$-coaction on $V_{S^{2N-1}}$ by $\DEL_{S^{N-1}}$, we see that if $(\DEL_{SU_N,\b} \oby \DEL_{S^{2N-1}})(\w) = \w \oby 1$, then clearly
\begin{align*}
\sum_{i \in \{+,0,-\}} (\DEL_{SU_N,\b} \oby \DEL_{S^{2N-1}})(w_i) ~ = \sum_{i \in \{+,0,-\}} w_i \oby 1.
\end{align*}
Thus, the lemma would follow if we could show that $\DEL_{S^{2N-1}}(V_i)  \sseq V_i \oby \bC_q[SU_{N-1}]$, for  all $i = +,0,-$.
That this is true is verifiable by direct calculation: For $i \neq 1$,
\begin{align*}
 \DEL(\ol{u^i_1}) & =  \sum_{k=1}^N \ol{u^k_1} \oby S(\b_N(u^i_k)) = \sum_{k=2}^N \ol{u^k_1} \oby S(u^{i-1}_{k-1}),
\end{align*}
and so, $\DEL_{S^{2N-1}}(V_+) \sseq V_+ \oby \bC_q[U_{N-1}]$. Moreover,
\begin{align*}
 \DEL(\ol{S(u^1_i)}) & =  \sum_{k=1}^N \ol{S(u^1_k)} \oby S(\b_N(S(u^k_i))) = \sum_{k=2}^N \ol{S(u^1_k)} \oby S^2(u^{k-1}_{i-1}),
\end{align*}
and so, $\DEL_{S^{2N-1}}(V_-) \sseq V_- \oby \bC_q[U_{N-1}]$. Finally, for $\bC e^0$,
\[
\DEL_{S^{2N-1}}(\ol{u^1_1 -1}) = \sum_k \ol{u^k_1} \oby S(\b_N(u^1_k)) - 1 \oby 1 = \ol{u^1_1} \oby 1 - 1 \oby 1 = e^0 \oby 1,
\]
and so, $\DEL_{S^{2N-1}}(\bC e^0) \sseq \bC e^0 \oby \bC_q[U_{N-1}]$.
\qed

As we shall see in the proof of the next proposition, $\bC e^0$ is not a right submodule, and so, $\Om^1_q(S^{2N-1})$ is not a right submodule. Thus, the above decomposition is {\em not} a decomposition into subcalculi.
\begin{prop} \label{prop:ecommrels}
The following right module relations hold in $\Om_q^1(S^{2N-1})$:
\begin{align*}
e^+_i z_r = q^{1-\frac{2}{N}}z_r e^+_i, &  &  e^+_iz^*_r = q^{\frac{2}{N}-1}z^*_r e^+_i, & &  e^-_iz_r = q^{1-\frac{2}{N}}z_r e^-_i, &  &  e^-_iz^*_r = q^{\frac{2}{N}-1} z^*_r e^-_i,
\end{align*}
\begin{equation*}
e^0z_r =  q^{2-\frac{2}{N}}z_re^0 + (q^{2-\frac{2}{N}}-1)\sum_{k=2}^{N} u^r_{k}e^+_{k-1},
\end{equation*}
\begin{equation*}
e^0z^*_r =  q^{\frac{2}{N}-2}z_r^*e^0 +  (q^{1+\frac{2}{N}})(q^{\frac{2}{N}}-q^2)  \sum_{k=2}^{N}  q^{-2k} S(u_r^{k}) e^-_{k-1}.
\end{equation*}
\end{prop}

\demo We shall only treat the actions of $z_r$ on $e^+_i$ and $e^0$, since the derivation of each of the other actions is directly analogous to one of these two. From (\ref{somegahbyh+}) we have that
\[
e^+_i z_r = (1 \oby \ol{u^i_1})u^r_1 = \sum_{k=1}^N u^r_k \oby \ol{u^i_1u^k_1}.
\]
Using (\ref{q1}) and (\ref{q2}), it is routine to show that $\ol{u^i_1u^1_1} = q^{(1-\frac{2}{N})}\ol{u^i_1}$, and  $\ol{u^i_1u^k_1} = 0$, for $k = 2, \ldots, N$. Thus, 
\[
e^+_iz_r =  q^{(1-\frac{2}{N})} u^r_1 \oby \ol{u^i_1} = q^{(1-\frac{2}{N})} z_r e^+_i.
\]
For the action of $z_r$ on $e^0$, we have that
\begin{align*}
e^0 z_r = (1 \oby (\ol{u^1_1-1}))u^r_1 = \sum_{k=1}^N u^r_k \oby \ol{(u^1_1 -1)u^k_1} = \sum_{k=1}^N u^r_k \oby \ol{u^1_1 u^k_1- u^k_1}.
\end{align*}
Using  (\ref{q1}) and (\ref{q2}) again, it is easy to show that $\ol{u^1_1u^1_1 - u^1_1} = q^{2-\frac{2}{N}}\ol{u^1_1 - 1}$,  and $\ol{u^1_1u^k_1 - u^k_1} = (q^{2-\frac{2}{N}} - 1)\ol{u^k_1}$, \, for $k = 2, \ldots, N$. Thus,
\begin{align*}
e^0.z_r ~ & = ~ q^{2-\frac{2}{N}}u^r_1\oby(\ol{u^1_1-1}) +  (q^{2-\frac{2}{N}}-1)\sum_{k=2}^{N} u^r_k\oby \ol{u^k_1} \\
          & = ~ q^{2-\frac{2}{N}}z_re^0 +  (q^{2-\frac{2}{N}}-1)\sum_{k=2}^{N} u^r_ke^+_k.
\end{align*}
 \qed

Finally, we establish the following lemma which gives an explicit description of the derivative for $\Om^1_q(\ccsn)$.

\begin{prop}\label{dforsphere}
In $\Om_q^1(S^{2N-1})$ it holds that
\begin{align}
 \exd z_i = z_i e^0 + \sum_{k=1}^{N-1} u^{i}_{k+1} e_{k}^+, &   &  \exd z_i^* = -q^{\frac{2}{N}-2}z^*_i e^0 - q^{1+\frac{4}{N}}\sum_{k=1}^{N-1} q^{-2(k+1)} S(u^{k+1}_i) e^-_{k}.
\end{align}
\end{prop}

\demo Calculating the formula for $\exd z_i$ is quite routine:
\[
\exd u^i_1  = \sum_{k=1}^{N} (u^i_k \oby \ol{u^k_1} - u^i_1 \oby 1)  = u^i_1 \oby \ol{(u^1_1 -1)} + \sum_{k=2}^{N} u^{i}_{k} \oby \ol{u^k_1} = u^i_1 e^0 + \sum_{k=2}^{N} u^{i}_{k} e_{k-1}^+.
\]
The derivation of the formula for $\exd z^*_i$ is slightly more involved: First we note that 
\begin{align*}
\exd (S(u_i^1)) =  \sum_{k=1}^{N} S(u^k_i) \oby \ol{S(u^1_k)} - S(u_i^1) \oby 1 =  S(u_i^1) \oby \ol{S(u^1_1) -1} + \sum_{k=2}^{N} S(u^k_i) \oby \ol{S(u^1_k)}.
\end{align*}
We then recall that $\ol{S(u^1_k)} = -q^{1+\frac{4}{N} -2k}\ol{u^1_k}$, and use (\ref{qS}) to derive the relation  \linebreak 
$\ol{S(u^1_1)-1} = -q^{\frac{2}{N}-2}(\ol{u^1_1 - 1})$. \qed

\section{A Framing for the Quantum Projective Spaces}\label{Sec:CPN}

Just as for the quantum spheres, we will now look at the framing of $\Om^1_q(\ccpn)$ given by Theorem \ref{homfra}. We begin by calculating the action of the soldering form $\ta$ on the basis elements of $V_{\ccpn}$ given in Lemma \ref{lem:cpnbasis}:
\begin{align} 
\ta(\ol{z_{i1}}) & = m \circ  (S \oby \exd) (\sum_{k,l=1}^N u^i_kS(u^l_1) \oby u^k_1S(u^1_l)) = \sum_{k,l=1}^N q^{2(l-1)}u^l_1S(u^i_k)\exd z_{kl}; \label{cpnsthetaofe1}\\
\ta(\ol{z_{1i}}) & = m \circ  (S \oby \exd) (\sum_{k,l=1}^N u^1_kS(u^l_i) \oby u^k_1 S(u^1_l)) = \sum_{k,l=1}^N q^{2(l-i)}u^l_iS(u^1_k)\exd z_{kl}. \label{cpnsthetaofe2}
\end{align}
Moreover, from (\ref{ztou}), we also have $\ta(\ol{z_{i1}}) = q^{\frac{2}{N}-1}e^+_i$, and $\ta(\ol{z_{1i}}) = - q^{3+\frac{2}{N}-2i}e^-_i$.

\subsubsection*{Holomorphic and Anti-Holomorphic Calculi} 

We shall now see that the framing presented above gives us a canonical decomposition of $\wcpn$ into a direct sum of two subcalculi. This decomposition  generalises the decomposition of the cotangent space of $\ccpn$ into its holomorphic and anti-holomorphic parts. We omit the proof
which is exactly analogous to the proof of Lemma \ref{decompositionofOm}.
\begin{prop} \label{cpndecomp}
Denoting 
\begin{align*}
\Om^{(1,0)}_q(\ccpn) = (\csun \oby V_+)^{U_{N-1}}, & &  \Om^{(0,1)}_q(\ccpn) = (\csun \oby V_-)^{U_{N-1}}, 
\end{align*}
we have a vector space decomposition:
\[
\Om^1_q({\bC P^{N-1}}) = \Om^{(1,0)}_q(\ccpn) \oplus \Om^{(0,1)}_q(\ccpn).
\]
\end{prop}

We recall that the decomposition of $\Om^1_q(S^{2N-1})$ into subspaces given in Lemma \ref{decompositionofOm} is not a decomposition into subcalculi. The above decomposition, however, {\em is}:

\begin{cor}
Let $\proj_+$ and $\proj_-$ be the canonical projections from $\Om^1_q(\bC P^{N-1})$ to $\Om^{(1,0)}_q(\ccpn)$ and $\Om^{(0,1)}_q(\ccpn)$ respectively,
and let us denote 
\begin{align*}
\del = \proj_+ \circ \exd,  & &   &  & \ol{\del} = \proj_- \circ \exd.
\end{align*}
Both pairs $(\Om^{(1,0)}_q(\ccpn),\del)$  and $(\Om^{(0,1)}_q(\ccpn),\ol{\del})$
are subcalculi of $\Om^1_q(\bC P^{N-1})$. We call them the {\em holomorphic calculus} and {\em anti-holomorphic calculus} respectively.
\end{cor}

\demo 
Proposition \ref{prop:ecommrels} implies that $\Omega_q^{(1,0)}(\ccpn)$ and $\Omega_q^{(0,1)}(\ccpn)$ are both right submodules. The direct sum decomposition of Proposition \ref{cpndecomp} implies that $\exd = \del + \ol{\del}$, and so, $\Omega_q^{(1,0)}(\ccpn)$ and $\Omega_q^{(0,1)}(\ccpn)$ are spanned by elements of the form $f\del g$ and $f\ol{\del} g$ respectively, for $f,g \in \cpn$. All that remains to show is that $\del$ and $\ol{\del}$ satisfy the Leibniz rule: From the Leibniz rule for $\exd$ we have that $\exd (fg)  = (\exd f)g +  f\exd g$. Since $\exd = \del + \ol{\del}$, this means that
\begin{align*}
\del(fg) + \ol{\del}(fg) = (\del  f)g + (\ol{\del} f)g + f\del g  + f(\ol{\del} g).
\end{align*}
But $\Om^{(1,0)}_q(\ccpn)$ and $\Om^{(0,1)}_q(\ccpn)$ are both right submodules, so our direct sum decomposition says we must have
\begin{align*}
\del(fg) = (\del f)g +  f\del g, &  & \ol{\del}(fg) = (\ol{\del} f)g +  f\ol{\del} g.
\end{align*}  \qed

We finish this section by finding formulae for the actions of $\del$ and $\ol{\del}$ on the elements of $\cpn$, as well as explicit descriptions of the right module relations of $\Om^{(1,0)}_q(\ccpn)$ and $\Om^{(0,1)}_q(\ccpn)$.
\begin{cor} \label{holoantiholocommrels}
The holomorphic and anti-holomorphic derivatives act on $\cpn$ according to:
 \begin{align}\label{delactions}
\del z_{ij} = q^{(\frac{2}{N}-1)} \sum_{k=2}^{N} u^i_{k}S(u^1_j) e^+_{k-1}, &  &  \ol{\del}z_{ij} = -q^{(3+\frac{2}{N})} \sum_{k=2}^{N} q^{-2k} u^i_1S(u^{k}_j) e^-_{k-1}.
 \end{align}
Moreover, we have the relations:
\begin{align} 
(\del z_{ij})z_{rs} & = \sum_{a,b,c,d,e,f=1}^N q^{\l_{bjs}}R^{ja}_{rb}\ol{R}^{ai}_{cd}R^{ec}_{fs} z_{de} \del z_{fb},
\label{thedelrelations}\\
(\ol{\del}z_{ij})z_{rs}   & = \sum_{a,b,c,d,e,f=1}^N q^{\ol{\l}_{bki}} R^{ra}_{jb}\,R^{sa}_{cd} \, \ol{R}^{ie}_{df} \, z_{be}\ol{\del} z_{fc}, \label{thedelbarrelations}
\end{align}
where $\l_{bjs} = 2(b-j)+\sgn(b-s)-1$, with $\sgn$ the sign function, and  $\ol{\l}_{bki}= 2(b-r) + \sgn(b-i) + 1$.
\end{cor}
\demo  We begin by noting that 
\begin{align*}
\exd z_{ij}  & =\exd(u^i_1S(u^1_j))   = \sum_{k,l=1}^{N} (u^i_kS(u^l_j) \oby \ol{u^k_1S(u^1_l)}) - u^i_1S(u^1_j) \oby 1\\
             & = u^i_1S(u^1_j) \oby \ol{u^1_1S(u^1_1)-1} + \sum_{k = 2}^{N} u^i_kS(u^1_j) \oby \ol{u^k_1S(u^1_1)} + \sum_{l=2}^{N} u^i_1S(u^l_j) \oby \ol{(u^1_1S(u^1_l)}\\
             & = \sum_{k =2}^{N} u^i_kS(u^1_j) \oby \ol{z_{k1}} + \sum_{l=2}^{N} u^i_1S(u^l_j) \oby \ol{z_{1l}}.
\end{align*}
As we saw earlier, $\ol{z_{k1}} = q^{(\frac{2}{N}-1)} e^+_{k-1}$ and $\ol{z_{1l}} = -q^{(3+\frac{2}{N}-2l)} e^-_{l-1}$, and so, we have that
\begin{align*}
\exd z_{ij} = q^{(\frac{2}{N}-1)} \sum_{k = 2}^{N} u^i_kS(u^1_j) \oby e^+_{k-1} - q^{(3+\frac{2}{N})}\sum_{l=2}^{N} q^{-2l}u^i_1S(u^l_j) \oby e^-_{l-1}.
\end{align*}
The formulae for the actions of the operators $\del$ and $\ol{\del}$ now follow directly.

Proposition \ref{prop:ecommrels} implies that the generators $z_{ij}$ commute with the $e^{\pm}_i$. This means that the right module relations can be determined using only the relations of $\csun$. 
We omit the calculations which are routine, if quite tedious.
\qed

\begin{eg}
Let us see  what the holomorphic and anti-holomorphic calculi look like for the case of $N=2$:
From Theorem \ref{homfra} we have that the coaction $\DEL_{\bC P^1}$ acts on $V_+$ and $V_-$ according to
\begin{align}
\DEL_{\bC P^1}(\ol{z_{21}}) = z_{21} \oby t^{-2}, & & \DEL_{\bC P^1}(\ol{z_{12}}) = z_{12} \oby t^2.
\end{align}
Thus, we must have $\Omega^{(1,0)}(\bC P^1) \simeq \E_{2}$, and $\Omega^{(0,1)}(\bC P^1) \simeq \E_{-2}$. From (\ref{sphererels}) we see that $\bC_q[\bC P^1]$ is generated by the three elements $z_{12},z_{21}$, and $z_{22}$. Using (\ref{delactions}) we can calculate that
\begin{align*}
\del z_{12}       & =  -q \inv b^2e^+_1,  &  \del z_{21}       & = d^2e^+_1,     & & \del z_{22}       = -q^{-2} bde^+_1,\\
\ol{\del} z_{12} & = -  a^2e^-_1,        &  \ol{\del}z_{21} & = qc^2e^-_1,    & & \ol{\del}z_{22}  = - q^{-1}ac e^-_1.
\end{align*}
Moreover, (\ref{thedelrelations}) and (\ref{thedelbarrelations}) give us the relations
\begin{align*}
\del z_{12}\begin{cases}z_{12}\\ z_{21}\\ z_{22}\end{cases}
& = \begin{cases}q^{-2}z_{12}{\del} z_{12}\\
q^{2} z_{21}{\del} z_{12}\\
z_{22}{\del} z_{12},\end{cases}
&
\del z_{21}\begin{cases}z_{12}\\ z_{21}\\
z_{22}\end{cases}
& =\begin{cases}q^{-2}z_{12}\del z_{21}
-\mu z_{21}\del z_{12}\\ q^{-2} z_{21}\del z_{12}\\
q^{-4}z_{22}\del z_{12},\end{cases} \\
  &  &   &                       \\
\ol{\del} z_{12}\begin{cases}z_{12}\\ z_{21}\\ z_{22}\end{cases}
& =\begin{cases}q^{2}z_{12}{\del} z_{12}\\
q^{2} z_{21}{\ol{\del}} z_{12} + \mu z_{12}{\ol{\del}} z_{21}\\
q^{4}z_{22}{\ol{\del}} z_{12},\end{cases}
& \ol{\del} z_{21}\begin{cases}z_{12}\\ z_{21}\\
z_{22}\end{cases}
& =\begin{cases}q^{-2}z_{12}\ol{\del} z_{21}\\ q^{2} z_{21}\ol{\del} z_{21}\\
z_{22}\ol{\del} z_{21}.\end{cases}
\end{align*}
where $\mu = q^2 - q^{-2}$. Similar relations hold for $\del z_{22}$ and $\ol{\del} z_{22}$.
This recovers (in our conventions) the description of the Podle\'s calculus given in \cite{Maj}.
\end{eg}

\subsubsection{Relationship with the Heckenberger--Kolb Calculus}

We will now show that $\Om^{(1,0)}_q(\ccpn)$ and $\Om^{(0,1)}_q(\ccpn)$ are none other than the two calculi identified by the classification result of Heckenberger and Kolb discussed in the introduction.

Before we can do this, however, we need to recall two definitions: Firstly, a left-covariant first-order calculus over an \alg $A$ is called {\em irreducible} if it does not possess any non-trivial quotients by a left-covariant $A$-bimodule.
Secondly, we recall from \cite{Herm} the definition for dimension of a calculus that is used in the statement of the classification: Let $G \iito M$ be a quantum homogeneous space, and $\Om^1(M)$ a first-order differential calculus over $M$ with corresponding sub-bimodule $N_M$. The {\em dimension} of $\Om^1(M)$ is defined to be $\dim(M^+/R_{\Om^1(M)})$, where
\[
R_{\Om^1(M)} = (\e_H \oby (\id - \e_H))(N_M).
\]
As one might expect, this definition coincides with our notion of dimension for the special case of a quantum principal homogeneous space. To see this we first note that
\begin{align*}
(\e_H \oby (\id - \e_H))(G N_M) = \e(G) R_{\Om^1(M)} = R_{\Om^1(M)}.
\end{align*}
Now if $I_M \sseq M^+$ is the right ideal corresponding to $\Om^1(M)$, then it can be shown that $s(GN_M) = G \oby I_M$ (see \cite{MMFPhD} for details). Moreover, for any $g \oby \ol{v} \in G \oby I_M$, it holds that 
\begin{align*}
(\e_H \oby (\id - \e_H))(s(g \oby \ol{v})) &  = (\e_H \oby (\id - \e_H))(g S(v\1) \oby \ol{v\2}) \\  
& =   \e(g) \e(v\1)  (v\2 - \e(v\2))\\ 
& = \e(g)v.
\end{align*}
Hence, we have 
\begin{align*}
(\e_H \oby (\id - \e_H))(G N_M) = (\e_H \oby (\id - \e_H))(s(G \oby I_M)) = \e(G)I_M = I_M.
\end{align*}
From this we can conclude that $R_{\Om^1(M)} = I_M$, and that the two notions of dimension do indeed coincide.

Let us now state the classification result for the special case of $\bC_q[\bC P^{N-1}]$:
\begin{thm} \cite{HK} \label{HKClass}
There exist exactly two non-isomorphic finite-dimensional irreducibe left-covariant first-order differential calculi over quantum projective $(N-1)$-space. Each has dimension $N-1$.
\end{thm}
Since both $\Om^{(1,0)}_q(\ccpn)$ and $\Om^{(0,1)}_q(\ccpn)$ have dimension $N-1$, they must both be irreducible (since otherwise there would exist an irreducible left-covariant calculus of dimension strictly less than $N-1$ in contradiction of the theorem). Moreover, it is easy to see that $\Om^{(1,0)}_q(\ccpn)$ and $\Om^{(0,1)}_q(\ccpn)$ correspond to different ideals of $\cpn^+$, and consequently are non-isomorphic. This gives us the following corollary:
\begin{cor}
The two calculi identified in Theorem \ref{HKClass} are $\Om^{(1,0)}_q(\ccpn)$ and \linebreak $\Om^{(0,1)}_q(\ccpn)$.
\end{cor}

\section{Connections}

In this section we shall discuss connections for the bundles $\a_N:\csun \to \cpn$ and $\csn \hookleftarrow \cpn$. Beginning with $\a_N:\csun \to \cpn$, we propose
\begin{align*}
i: \Lambda^1_{U_{N-1}} \to \Lambda^1_{SU_N}, &  &\ol{\dt_{N-1} \inv - 1} \mto e^0,
\end{align*}
as a bicovariant splitting map. That $i$ satisfies $\a_N \circ i = \id$ is obvious, while (\ref{bicovcond}) follows from 
\begin{align*}
\ol{\Ad_{R,SU_N}} \circ i (\ol{\dt_N \inv -1}) & = \ol{\Ad_{R,SU_N}}(e^0) = \sum_{k,l=1}^N (\ol{u^k_l} \oby \a_N(S(u^1_k)u^l_1)) - 1 \oby 1 \\ 
 & = \ol{u^1_1} \oby 1                                                       - 1 \oby 1  = e^0 \oby 1  = i(\ol{\dt_N \inv -1}) \oby 1 \\ & = (i \oby \id) \circ \ol{\Ad_{R,U_{N-1}}}(\ol{\dt_N \inv -1}).
\end{align*}
We denote the connection form corresponding to $i$ by $\w$. The following result establishes some properties of the connection $\Pi_\w$ associated to $i$:
\begin{lem} \label{sunconnectionformaction}
The connection $\Pi_\w:\Om^1_q(SU_N) \to \Om^1_q(SU_N)$ is strong and satisfies 
\begin{align*}
\Pi_\w(e^0) = e^0, & & {\Pi_\w(e^+_i) = \Pi_\w(e^-_i) = 0}, & & (i = 1, \ldots, N-1).
\end{align*}
\end{lem}

\demo From (\ref{mver}) and (\ref{pifromom}) we see that
\begin{align*}
\Pi_\w(e^0) = m \circ (\id \oby \w) \circ \ver (\ol{u^1_1-1}) =  \w \circ (\ol{u^1_1-1}) = \w(\ol{\dt_{N-1}\inv -1})=e^0.
\end{align*}
Similarly, $\a_N(u^{i+1}_1) = \a_N(u^1_{i+1}) = 0$ implies that $\Pi_\w(e^+_i) = \Pi_\w(e^-_i) = 0$. 

To show that $\Pi_\w$ is strong, we must establish that 
\[
(\id - \Pi_\w)(\exd (\csun)) \sseq \csun \Om^1_q(\ccpn).
\]
From our results for the action of $\Pi_\w$, it is clear that $(\id - \Pi_\w)$ acts on a general form $\sum_{i=1}^{N-1} f^+_ie^+_i + f^0e^0 + \sum_{i=1}^{N-1} f^-_ie^-_i$, for $f^+_i,f^0,f^-_i \in \csun$, to give $\sum_i f^+_i e^+_i + \sum_i f^-_i e^-_i$. Thus, if we could show that ${e^{\pm}_i \in \csun \Om^1_q(\ccpn)}$, then it would follow that $\Pi_\w$ was strong. But this is a direct consequence of (\ref{cpnsthetaofe1}) and (\ref{cpnsthetaofe2}), and so, $\Pi_\w$ is indeed strong. \qed

A natural question to ask is whether or not $\Pi_\w$ restricts to a connection for the bundle $\csn \hookleftarrow \cpn$. The following lemma shows that it does.
\begin{prop}
The connection $\Pi_\w$ restricts to a connection for the bundle \linebreak $\csn \hookleftarrow \cpn$.  Its corresponding connection form is the mapping 
\begin{align*}
\w':\Lambda^1_{U_1} \to \Omega^1_q(S^{2N-1}), & & & \ol{t \inv - 1} \mto e^0.
\end{align*}
\end{prop}  
\demo 
We denote by $J:\Lambda_{U_1}^1  \to \Lambda_{U_{N-1}}^1$ the unique linear map for which for $J(\ol{t\inv-1}) =  e^0$. It is easy to see from (\ref{cpninsphere}) that we have the following commutative diagram:
\begin{align*} 
\xymatrix{
\Omega_q^1(S^{2N-1}) \ar[d]_{\iota} \ar[rr]^{\ol{\ver}~~~} & &  \csn  \oby \Lambda_{U_1}^1   \ar[d]^{\iot \oby J}\\
\Omega^1_q(SU_N)                    \ar[rr]_{\ol{\ver}~~~} & &  \csun \oby \Lambda^1_{U_{N-1}}.
}
\end{align*}
Moreover, the diagram 
\begin{align} \label{restconnform}
\xymatrix{ \Lambda_{U_1}^1      \ar[d]_{J} \ar[rr]^{\w'}  & &  \Omega_q(S^{2N-1}) \ar[d]^{\iota}\\
           \Lambda_{U_{N-1}}^1                  \ar[rr]_{\w}   & &  \Omega_q(SU_N),}
\end{align}
is also commutative, which means that we  have
\[
\Pi_{\w}|_{\Omega(\ccsn)} = m  \circ (\id \oby \w') \circ \ol{\ver}.
\]
Thus, if we could show that $\w'$ was a connection form, the lemma would follow. Using a commutative diagram argument this can be deduced from the fact that $\w$ is a connection form. However, it is much more economical to verify the requirements directly: That $\ver \circ \w' = 1 \oby \id$ follows from 
\begin{align*}
\ver \circ \w'(\ol{t\inv -1})          & = \ver (e^0) = 1 \oby \a_N(u^1_1-1) = 1 \oby \ol{t\inv-1};
\end{align*}
while $(\w' \oby \id)\circ \ol{\Ad_{R,U_1}} = \DEL_{\ccsn, \g}\circ \w'$ follows from the calculation
\begin{align*}
(\w' \oby \id)\circ \ol{\Ad_{R,U_1}}(\ol{t\inv-1}) & = (\w'(\ol{t\inv-1})) \oby 1 = e^0 \oby 1,
\end{align*}
and the fact that
\begin{align*}
\DEL_{\ccsn, \g} \circ \w'(\ol{t\inv-1}) & = \DEL_{\ccsn, \g}(e^0) = s(\sum_{k,l=1}^N 1 \oby \ol{u^k_l} \oby \g_N(S(u^1_k)u^l_1) - 1 \oby  1 \oby 1)\\
         & = s(1 \oby \ol{u^1_1 -1} \oby 1) = e^0 \oby 1.
\end{align*}
\qed

\begin{prop} \label{prop:associatedun-1}
Every associated vector bundle to $\csn \hookleftarrow \cpn$ is also an associated vector bundle to $\a_N:\csun \to \cpn$.
\end{prop}

\demo For any $\bC[U_1]$-coaction $\DEL_{V,U_1}$ on a vector space $V$,  we can define a $\bC_q[U_{N-1}]$-coaction on $V$ by setting $\DEL_{V,U_{N-1}} = (\id \oby j) \circ \DEL_{V,U_1}$. 
We will establish the proposition by showing that
$$
(\csn \oby V)^{U_1} = (\csun \oby V)^{U_{N-1}}.
$$
We begin by noting that (\ref{cpninsphere}) implies that the following diagram is commutative:
\begin{displaymath}
\xymatrix{\csn  \oby V    \ar[rrrr]^{\DEL_{S^{2N-1},\g} \oby \DEL_{V,U_1}} \ar[d]_{\iota }& & & &  \csn  \oby V \oby \bC_q[U_{1}] \ar[d]^{\iota \oby \id \oby j} \\
          \csun \oby V    \ar[rrrr]_{\DEL_{SU_N,\a} \oby \DEL_{V,U_{N-1}}}               & & & & \csun \oby V \oby \bC_q[U_{N-1}].}
\end{displaymath}
From this we see that $(\csn \oby V)^{U_{1}} \sseq (\csun \oby V)^{U_{N-1}}$. It remains to establish the opposite inclusion. If $\sum_{i=1}^M f_i \oby v_i$ is an element of $(\csun \oby V)^{U_{N-1}}$, then 
\begin{align*}
\sum_{i=1}^M (f_i)_{(1)} \oby (v_i)_{(0)} \oby \a_N((f_i)_{(2)})j((v_i)_{(1)}) = \sum_{i=1}^M f_i \oby v_i \oby 1,
\end{align*}
where we denote $\DEL_{V,U_1}(v) = v_{(0)} \oby v_{(1)}$, for $v \in V$. Operating on both sides by ${\id \oby \id \oby \d_{N-1}}$ (where, as in Lemma \ref{spherelcoaction}, $\d_{N-1}$ is the canonical projection from $\bC_q[U_{N-1}]$ to $\bC_q[SU_{N-1}]$) gives 
\begin{align*}
\sum_{i=1}^M (f_i)_{(1)} \oby (v_i)_{(0)} \oby (\d_{N-1} \circ \a_N((f_i)_{(2)}))(\d_{N-1} \circ j((v_i)_{(1)}) = \sum_{i=1}^M f_i \oby v_i \oby 1.  
\end{align*}
Recalling that $\d_{N-1} \circ \a_N = \b_N$, and noting that $\d_{N-1} \circ j(f) = \e(f)1$, for all $f \in \bC[U_1]$, we see that 
\begin{align*}
\sum_{i=1}^M (f_i)_{(1)} \oby v_i \oby \b_N((f_i)_{(2)}) = \sum_{i=1}^M f_i \oby v_i \oby 1.
\end{align*}
Assuming, without loss of generality, that the $v_i$ are linearly independent, we get
\begin{align*}
(f_i)_{(1)} \oby \b_N((f_i)_{(2)}) = f_i \oby 1,&  &  (i =  1, \ldots, M).
\end{align*}
Since $(f_i)_{(1)} \oby \b_N((f_i)_{(2)}) = \DEL_{SU_N,\b}(f_i)$, we must have $f_i \in \csn$, for all $i$.

\qed

\begin{cor}
The connection $\Pi_{\w'}$ associated to $\w'$ is strong. Moreover, the action of the covaraint derivative $\nabla_{\w'}$ induced by $\w'$ coincides with the action of $\nabla_{\w}$ on all the associated vector bundles of $\csn \hookleftarrow \cpn$.
\end{cor}

\demo
Recall that we have a $\bZ$-grading on $\csn$ induced by the $\bC[U_1]$-coaction $\DEL_{\ccsn,\g}$. As before, we denote the space of homogeneous elements of degree $p$ by $\E_p$. Since $\E_p$ is an associated bundle to $\csn \hookleftarrow \cpn$, the previous proposition tells us that it is also an associated bundle to $\a_N:\csun \to \cpn$. Since $\w$ is strong, we have $\nabla_\w(\E_p) \sseq \E_p \oby \Om^1_q(\ccpn)$, for all $p$, and consequently
\[
(\id - \Pi_{\w}) \circ \exd (\bC_q[\ccsn]) \sseq \bC_q[\ccsn]\Om^1_q(\ccpn).
\]
Now $\id - \Pi_{\w'}$ is equal to the restriction of  $\id - \Pi_{\w}$ to $\Om^1_q(\ccsn)$, and so, we must have that $\w'$ is strong.

That the actions of the covaraint derivatives $\nabla_{\w}$ and $\nabla_{\w'}$ coincide is obvious from (\ref{covderivative}).

\qed

\begin{eg}

As an easy example let us look at the action of $\nabla$ on the simplest element of the quantum line  bundle $\E_{1}$. For $z^*_i \in \E_{1}$, we recall that
\[
\exd z_i^* = -q^{\frac{2}{N}-2}z^*_i e^0 - q^{1+\frac{4}{N}}\sum_{k=2}^{N} q^{-2k} S(u^{k}_i) e^-_{k-1}.
\]
Applying $\Pi-\id$ gives us
\begin{align*}
 \nabla(z^*_i) & = - q^{1+\frac{4}{N}}\sum_{k=2}^{N} q^{-2k} S(u^{k}_i) e^-_{k-1} = - q^{1+\frac{4}{N}}\sum_{k=2}^{N}\sum_{l=1}^N  q^{-2k} S(u^1_l)u^l_1 S(u^{k}_i) e^-_{k-1}\\
               & = q^{(1+\frac{4}{N})-(3+\frac{2}{N})}\sum_{l=1}^N (S(u^1_l)(-q^{3+\frac{2}{N}})(\sum_{k=2}^{N} q^{-2k}u^l_1
                S(u^{k}_i) e^-_{k-1})) \\
               & = q^{\frac{2}{N}-2}\sum_{l=1}^N z^*_l \oby \ol{\del} z_{li}.
\end{align*} 
\end{eg}

\bigskip
School of Mathematical Sciences, Queen Mary, University of London,
327 Mile End Rd, London E1 4NS, England

{\em e-mail}: \tt{reob@maths.qmul.ac.uk}

\end{document}